July14, 2015
v2: September 4, 2015
v3: November 4, 2015
v4: November 20, 2015
Yu-Chiu Chao
TRIUMF
Vancouver, Canada

# A Deterministic Program for Obtaining Optima under Constraints for Any Analytical System

## Abstract

Conceptual framework is laid out of a deterministic program capable of obtaining optimum solutions with or without constraints for any reasonably behaved analytical system. Recipe implementable as a well-behaved Runge-Kutta procedure is given. Determinism means no inspired initial guesses or random number trials. The program follows well-defined steps between well-defined start and end points, to a large extent configuration independent. This program is also conjectured to lead to global optimum solutions based on evidence, short of a proof. Application to realistic problems is given as example.

## Introduction

The celebrated recipe of Lagrange [1] has been for 2 centuries canon and staple to people working on optimization under constraints in an analytical system. It is found in all entry-level calculus textbooks and requires no introduction beyond the basic formulation:

$$\begin{cases} \nabla F = \lambda \cdot \nabla H \\ H = h \end{cases} \rightarrow \; k_1^O, k_2^O, k_3^O, \dots k_N^O, \lambda^O \rightarrow \; F = f(h) \tag{1}$$

In the above we have a system objective to be optimized, denoted by a function $F$ of variables $k_m$, subject to a constraint function $H = h$.

Variables: $k_1, k_2, k_3, \dots k_N$

Objective function: $F(k_1, k_2, k_3, \dots k_N)$ to be optimized

Constraint: $H(k_1, k_2, k_3, \dots k_N) = h$

The solution is obtained by solving (1) with $N+1$ equations and $N+1$ unknowns, $k_m$ and $\lambda$, with optimal solution for $F$ denoted $f$ in (1). This prescription does not specify the algorithm by which the solution is realized. It is up to the user to choose such an algorithm, as well as accept all its caveats and pitfalls, which can be myriad as $F$ and $H$ grow in complexity, and the inherent uncertainties in these algorithms start to become an issue. In some cases artificial choice of an initial guess will determine the quality of the final answer, while in others the algorithm invokes random number search in hope of landing in the right neighborhood of solution space.

In a slightly modified view of the problem, one may ask instead "As $h$ is varied, how does $f$ change in response while always satisfying (1)?" This amounts to mapping out the optimal path of trade-off between objective and constraint for not only one, but an entire range of such options. At every point on this path the optimal conditions for both $f$ and $h$ are simultaneously met. In other words, for any given $f$ on the curve, $h$ is at a local optimum, and vice versa. The distinction between objective and constraint in this view has disappeared, and $f$ and $h$ are simply two competing objectives on the same footing. One may therefore write the equivalent problem, trading objective for constraint:

$$\begin{cases} \nabla F = \lambda \cdot \nabla H \\ F = f \end{cases} \rightarrow \; k_1^O, k_2^O, k_3^O, \dots k_N^O, \lambda^O \rightarrow \; H = h(f) \tag{2}$$

We are actually interested in $k_m^O(h), \lambda^O(h)$ as functions of $h$ (or f: $k_m^O(f), \lambda^O(f)$). As $h$ varies, the solution traces out a 1D curve in the space spanned by $k_1, k_1, k_1, \dots k_N$. One can parametrize this curve with $h$, or $f$, or even $\lambda$. It is also clear from the above that as long as one stays on this 1D solution curve, any one of the trio ($h$, $f$, $\lambda$) uniquely[1] determines the other two.

The above observation can be taken further, namely, one can map out the 1D curve by varying $\lambda$, thereby reducing the numerical difficulty of the problem because one now has one less equation to satisfy[2].

$$\nabla F = \lambda \cdot \nabla H \; \rightarrow \quad k_1^O, k_2^O, k_3^O, \dots k_N^O \; \rightarrow \; \begin{cases} F = f(\lambda) \\ H = h(\lambda) \end{cases} \tag{3}$$

Here objective and constraint are on manifestly symmetric footing, with solutions $k_m^O(\lambda)$ (but no $\lambda^O(\lambda)$ now). This is numerically easier than the other formulations, especially when dimensionality and complexity increase in $F$ and $H$. In practice (3) is advantageous, and even essential, also because (a), one does not need a priori knowledge of where it should start and end in a certain class of problems, and (b), it is the only way to extricate the process out of singular points, as will be seen.

## *Solution Curve in k-Space, Symmetry between h and f, and Constraint Independent Objective Optimum*

The 1D solution curve in the space spanned by $k_1, k_1, k_1, \dots k_N$ can be parametrized by any one of the inter-dependent variables ($h$, $f$, $\lambda$). Formulations (1), (2), and (3) nominally describe the same curve. One can vary the variable in red in each formulation to map out this curve, although there is indispensability and certain numerical advantage in formulation (3) as noted in the footnote below.

The relations between $h$, $f$, and $\lambda$, assuming one stays on the 1D solution curve, deserve some delineation. It is obvious from (1) or (2) that

$$\left. \frac{df}{dh} \right| = \lambda \tag{4}$$

where the vertical bar indicates the derivative is taken along the 1D curve constrained by (1), (2), or (3) with the only DOF given by the variable in red in respective systems. The ratio between rates of change in $f$ and $h$ is thus given by $\lambda$. More specifically as one moves along the 1D solution curve $\lambda$ can change from one sign to the other, signaling growing or diminishing $f$ with increasing $h$. The point $\lambda = 0$ thus represents the point where $f$ reaches an extremum in a more absolute and constraint-independent[3] sense. This is the familiar absolute extremum for $f$, itself often the aim of more ambitious optimization programs[4]:

$$\nabla F = 0 \tag{5}$$

[1] For now we ignore cases of possible multiple solution paths corresponding to $F$ and $H$ tangent to each other at multiple locations in the $k$-space. Proof of uniqueness is beyond the current scope.

[2] In addition (1) and (2) can be numerically very unstable near $|k| = 0$. In general mapping out global behavior of both (1) and (2) suffers from a lack of well-defined staring or ending point, whereas for (3) these points correspond to $-\infty$ and 0.

In the objective-constraint context it has the remarkable property of being independent of the form of the constraint $H$ (A point is tangent to any curve or surface). This statement is less trivial than it sounds if one realizes that (5) corresponds to the point in the space of $k_m$, up to which relaxation of the constraint resulted in progressively improved objective, and beyond which the same relaxation causes the objective to degrade. This property holds for <u>all forms</u> of constraints. We will call this point the Constraint Independent Optimum (CIO).

As the boundary between objective and constraint is blurred, we will drop this nomenclature and in the following refer to $F$ and $H$ simply as two competing objectives, with respective CIO's given by (5) and (6) below.

$$\nabla H = 0 \tag{6}$$

As will be seen below, in some cases, this implies that a problem can be made more deterministically solvable by <u>adding</u> artificial constraints!

## *Relations between h, f, & $\lambda$ along the Solution Curve and Evolution of $k_m$*

Optimization based on (1) is most typically encountered. Its conjugate formulation (2) is not as frequently explored. Here we conceive a program incorporating (3) since we are interested in neither optimal $f$ for a given $h$ nor vice versa, but the evolution of simultaneously optimized $f$ and $h$ over the entire range of tradeoff. Solving (3) over the entire range of $\lambda$ is equivalent to accomplishing the following:

- Mapping out global tradeoff between $f$ and $h$ constrained by the optimum condition[5]
- Mapping out paths of $k_1, k_1, k_1, \dots k_N$ along the above curve.

In addition (3) affords two algorithmic and numerical advantages:

- (3) has one less equation to satisfy, which can be critical when degrees of freedom gets large,
- $\lambda$ has well-defined and universal starting and ending values over the full range of interest for a special class of problems, as will be seen. The same is not true with either $f$ or $h$,
- (3) provides a way to negotiate out of singularity traps of (1) or (2), which the latter, even taken together, cannot achieve.

The second point above can in principle be exploited, together with (5) and (6), to make it true for <u>all</u> classes of problems. All three points above point to the advantage of (3), the $\lambda$-based solution path, as basis to a deterministic numerical program.

To go about mapping out the evolution of (3), namely, the complete set of $k_1^O(\lambda), k_2^O(\lambda), k_3^O(\lambda), \dots k_N^O(\lambda)$ for each $\lambda$ between starting and ending values of interest, it would be useful to further explore the inter-dependency between $h$, $f$, and $\lambda$, along the 1D solution path in the spirit of (4). We note that along this path

---

[3] Progressive relaxation of constraint up to this point has led to continued improvement on the objective. Further relaxation beyond this point will result in diminished objective. This is true <u>regardless of the detail of the constraint of interest</u>.
[4] We will not consider the saddle-point possibility for (4) or (5) in this note, as despite being interesting in its own right, it distracts from the main purpose of the current discussion.
[5] This tradeoff is akin to the Pareto front in genetic optimization problems.

$$G(k, \lambda) = \nabla(F - \lambda \cdot H) = 0 \tag{7}$$

where $G(\lambda)$ is simply a vector function as defined. It is trivial that along this curve the following is also true[6]

$$\left.\frac{dG_i(k,\lambda)}{d\lambda}\right| = \sum_{j=1}^{N} \frac{dk_j(\lambda)}{d\lambda}\frac{\partial G_i}{\partial k_j} + \frac{\partial G_i}{\partial \lambda}$$
$$= \sum_{j=1}^{N} \frac{dk_j(\lambda)}{d\lambda}\frac{\partial^2(F - \lambda \cdot H)}{\partial k_i \partial k_j} - \frac{\partial H}{\partial k_i} = 0 \tag{8}$$

or

$$\left.\frac{d\boldsymbol{k}}{d\lambda}\right| = \boldsymbol{M}^{-1} \cdot \boldsymbol{R}$$
$$\boldsymbol{k} = \left(k_1^O(\lambda), k_2^O(\lambda), \dots k_N^O(\lambda)\right), \qquad \boldsymbol{M}_{ij} = \frac{\partial^2(F(\boldsymbol{k}) - \lambda \cdot H(\boldsymbol{k}))}{\partial k_i \partial k_j}, \qquad \boldsymbol{R}_i = \frac{\partial H(\boldsymbol{k})}{\partial k_i} \tag{9}$$

***M*** is known as the Hessian. We have another trivial relation along this path from, for example, (1)

$$\sum_{j=1}^{N} \frac{dk_j(\lambda)}{dh}\frac{\partial H}{\partial k_j} = \left.\frac{d\boldsymbol{k}}{dh}\right| \cdot \boldsymbol{R} = 1 \tag{10}$$

Substituting (9) into (10) and connecting $d\boldsymbol{k}/d\lambda$ and $d\boldsymbol{k}/dh$ through the scalar $d\lambda/d\,h$, we get

$$\left.\frac{d\lambda}{dh}\right| \cdot (\boldsymbol{R}^T \cdot \boldsymbol{M}^{-1} \cdot \boldsymbol{R}) = 1,$$
$$\left.\frac{d^2 f}{dh^2}\right| = \left.\frac{d\lambda}{dh}\right| = \left.\frac{1}{\boldsymbol{R}^T \cdot \boldsymbol{M}^{-1} \cdot \boldsymbol{R}}\right| \tag{11}$$

Equivalently one can use instead the constraint in (2) and arrive at a conjugate set of equations describing the same family of optimized solutions ($\mu = 1/\lambda$):

$$\left.\frac{d\boldsymbol{k}}{d\mu}\right| = \boldsymbol{N}^{-1} \cdot \boldsymbol{S}$$
$$\boldsymbol{k} = \left(k_1^O(\lambda), k_2^O(\lambda), \dots k_N^O(\lambda)\right), \qquad \boldsymbol{N}_{ij} = \frac{\partial^2\big(H(\boldsymbol{k}) - \mu \cdot F(\boldsymbol{k})\big)}{\partial k_i \partial k_j}, \qquad \boldsymbol{S}_i = \frac{\partial F(\boldsymbol{k})}{\partial k_i} \tag{12}$$

as well as

$$\left.\frac{dh}{df}\right| = \mu = \frac{1}{\lambda}$$
$$\left.\frac{d^2 h}{df^2}\right| = \left.\frac{d\mu}{df}\right| = \left.\frac{1}{\boldsymbol{S}^T \cdot \boldsymbol{N}^{-1} \cdot \boldsymbol{S}}\right| \tag{13}$$

This conjugate formulation proves to be important for the robustness of algorithm implementation as will be seen below (for example see (20)).

---

[6] Most of the following discussion can be summarized by the tensor relation $\nabla_{(k,\lambda)}\boldsymbol{G}| = 0$, with $\nabla_{(k,\lambda)}$ the vector differential operator acting on both $k_m$'s and $\lambda$. At some point for ease of notation the explicit dependence on $\lambda$ and $k$ is dropped from the formulas.

Equation (9) suggests that, with a suitably chosen initial set of $k_1, k_1, k_1, \dots k_N$ and $\lambda$, it can be integrated to map out the entire family of solutions for optimal $k_1, k_1, k_1, \dots k_N$ as functions of $\lambda$, and thus the global evolution of $f$ and $h$. In reality this can be hindered by numerical complexity of the right hand side of (9), and presence of inflections and extrema in these functions, as will be analyzed below.

Equation (11) or (13) are key to further defining the interplay between $h$, $f$, and $\lambda$. Staying with the formulation of (11), the derivatives between these quantities along the 1D solution curve is given in Table 1. These relations offer some useful insight on the process of integrating a system given by (9), at least conceptually, where one can envision starting from an initial condition, $k_1, k_1, k_1, \dots k_N$ already satisfying (3) at an initial $\lambda$, and integrating (9) to map out the entire family of solutions as $\lambda$ is varied over the range of interest. This solution path represents optimal[7] $f$ for any given $h$, and vice versa, resulting in a trade-off curve between $f$ and $h$, as shown in Figure 1. Relaxation on $f$ ($h$) leads to ever more optimal $h$ ($f$), until a point where further relaxation on $f$ ($h$) leads to a reversal in the trend for $h$ ($f$). This point of reversal for $h$ ($f$) corresponds to $\mu = 0$ ($\lambda = 0$) given by the CIO's (6) or (5).

By the above argument the range of interest for $\lambda$ is between $\pm\infty$[8] and 0, outside which the trade-off between optimal $f$ and $h$ ceases to exist and the reverse tradeoff takes place in which the relaxation of one drives the other further away from optimum[9]. This produces a well-defined range for $\lambda$ in system (3) independent of the details of the problem under investigation, as opposed to systems (1) or (2) where the exact range of the free parameter $f$ or $h$ cannot be known a priori, making it difficult to either find the initial condition for integrating (9), or know where the integration has reached the point of diminished return since the CIO has been reached.

## *Inflection and Extremum in the f–h Tradeoff and Alternating Integration for $k_m$*

In practice if (3) is followed, the integration of (9) will be stopped at local extrema of $\lambda$ before it reaches the endpoint. By (11) local minima in $\lambda$, signified by $d\lambda/d\,h = 0$, should never happen unless

$$Det(\boldsymbol{M}) = 0 \tag{14}$$

This condition is important in marking all "break points" in the process of integrating (9). $\lambda$ cannot be monotonically extended beyond such a point, so in order to continue with the integration, a different variable, such as $f$ based on (2) or $h$ based on (1), must be used subsequent to the $\lambda$–based procedure.

The point $Det(\boldsymbol{M}) = 0$ is itself interesting in that it corresponds to inflection points in the trade-off picture between $f$ and $h$, and thus can imply "diminished return" or "enhanced return" as explained in Figure 2. These points can provide insight on choosing a solution along the trade-off curve.

At the inflection point the integration has to continue on variables other than $\lambda$. One can use either $f$ or $h$[10]:

---

[7] For ease of discussion we assume the problem is to minimize $f$ for given $h$, and vice versa. In this regime $\lambda$ & $\mu$ are both negative.

[8] The sign is determined by the direction of tradeoff between $f$ and $h$.

[9] Of course the so called optimum is very much a human choice. One person's maximum is another's minimum. Everywhere along the solution curve, within and without this range, the Lagrangian condition holds pertaining to the extremum requirement.

[10] This stems from the fact, for example, that if one follows the system (1) with $f$ being the independent variable parametrizing the 1D solution curve, then $\lambda$ is a function of $f$, $\lambda^O(f)$, along the curve.

$$\left.\frac{d\boldsymbol{k}}{df}\right| = \frac{1}{\lambda}\cdot\frac{\boldsymbol{M}^{-1}\cdot\boldsymbol{R}}{\boldsymbol{R}^{\boldsymbol{T}}\cdot\boldsymbol{M}^{-1}\cdot\boldsymbol{R}} = \frac{1}{\lambda}\cdot\frac{Adj(\boldsymbol{M})\cdot\boldsymbol{R}}{\boldsymbol{R}^{\boldsymbol{T}}\cdot Adj(\boldsymbol{M})\cdot\boldsymbol{R}}$$
$$\left.\frac{d\boldsymbol{k}}{dh}\right| = \frac{\boldsymbol{M}^{-1}\cdot\boldsymbol{R}}{\boldsymbol{R}^{\boldsymbol{T}}\cdot\boldsymbol{M}^{-1}\cdot\boldsymbol{R}} = \frac{Adj(\boldsymbol{M})\cdot\boldsymbol{R}}{\boldsymbol{R}^{\boldsymbol{T}}\cdot Adj(\boldsymbol{M})\cdot\boldsymbol{R}} \tag{15}$$

where $Adj(\boldsymbol{M})$, the adjugate of $\boldsymbol{M}$, or transpose of the cofactor matrix of $\boldsymbol{M}$, is used,

$$Adj(\boldsymbol{M}) = Cof(\boldsymbol{M})^{\boldsymbol{T}} = Det(\boldsymbol{M})\cdot\boldsymbol{M}^{-1} \tag{16}$$

This makes (15) well behaved when (9) isn't, so that it can pick up after (9) when the latter is stopped due to $Det(\boldsymbol{M}) = 0$[11].

By the same token the integration along either $f$ or $h$ cannot be monotonically extended beyond extrema in $f$ or $h$. These correspond to points satisfying

$$\boldsymbol{R}^{\boldsymbol{T}}\cdot\boldsymbol{M}^{-1}\cdot\boldsymbol{R} = 0 \quad \text{or}$$
$$\boldsymbol{R}^{\boldsymbol{T}}\cdot Adj(\boldsymbol{M})\cdot\boldsymbol{R} = 0 \tag{17}$$

Note (17) does not necessarily make the numerator of (15) zero, as it does not require a rank-deficient $\boldsymbol{M}^{-1}$, but rather an $\boldsymbol{M}^{-1}$ whose eigenvalues are not all of the same sign[12].

Equations (9) and (15) indicate another significance of (3) as the principal means for mapping out the solution path when alternated with (1) or (2), rather than alternating between (1) and (2). This can be seen because (1) and (3), or (2) and (3) do not share the same set of singular points, (14) and (17). Thus alternating between (1) and (3), or between (2) and (3), ensures that the process does not get "stuck" at a singular point from which it cannot extricate itself. If on the other hand one uses (1) and (2), the process can get stuck when (17) happens and the path following either (1) or (2) cannot advance further. This is demonstrated by the solution path plotted in the $f$–$\lambda$ or $h$–$\lambda$ plane, showing rounded corners around these singular points that can be negotiated by alternating (3) with (1) or (2). The same points in the $f$–$h$ plane sit at cusps of intersecting curves with incompatible slopes[13], making smooth negotiation impossible. Figure 3 gives a graphical demonstration of this.

In view of the highly nonlinear RHS of (9) and (15) potentially making numerical integration difficult, one can instead successively solve (3) at small incremental values of $\lambda$ using preceding step as input. When an inflection point is encountered, switch to (1) or (2) at small incremental values of $f$ or $h$, until extrema is reached, where the process switches back to (3). This process is iterated until $\lambda = 0$ is reached.

The configuration insensitive nature of this procedure makes it well suited for large scale, mass production applications, such as generating interpolation database on widely varying system parameters. In some sense the procedure does not have free parameters to tweak, fine tune, or even guess in response to changing configuration parameters, nor is there a need to know a priori where the process should start or end as a function of problem configuration.

In summary, the proposed solution program amounts to alternating integration of (9) and (15), with switch-over points defined by inflection (14) or extrema (17), within range of $\lambda$ delimited by universal end points ±∞ and 0.

---

[11] Likewise (9) is well behaved when (15) is stopped due to (17) since (17) does not require rank deficiency.
[12] Actually, being the inverse of $\boldsymbol{M}$, $\boldsymbol{M}^{-1}$ cannot be rank-deficient.
[13] Since the slope in this case is just $\lambda$, it is always negative and can never "round the corner" by becoming 0 or positive.

## *Imposing Artificial Constraints to Make an Intractable Problem Solvable*

For the special class of constraints $H$, notably those involving only up to 2nd order polynomials in $k_m$ such as a quadratic sum of all $k_m$'s, (6) is algebraically easy and thus provides a robust starting point for the program. From this point on one only needs to take $\lambda$ from $\pm\infty$[14] to 0 by (3), including requisite alternation with (1) or (2), without concern over where the process should stop. The entire process is deterministic with universal start and end points: $\pm\infty$ and 0.

In some applications all one wants is the CIO without constraints, but the problem is intractable as numerical methods fail to produce a confident CIO. In such cases we can take advantage of the determinism of the constrained solution path and impose an artificial constraint $H$ with well-behaved and easily solvable CIO by (6), and integrate along the path to reach the true CIO at the other end for $F$[15].

To take the above argument even further, one can map out the solution path for an otherwise intractable constrained problem whose CIO is intractable on both ends and thus cannot provide either a suitable starting point or an ending point. This is done by first imposing an artificial constraint $H'$ with well-behaved and easily solvable CIO, integrating to the other end of the path to get the CIO for the true objective $F$, and then using this CIO as starting point to integrate backward along the true solution path under the desired constraint $H$ to the CIO on the other end. Figure 4 shows the reversibility of the integration process. The same tradeoff curve is mapped out from CIO on either end.

The above argument indicates that, conceptually at least, the procedure should deterministically solve arbitrary constrained optimization problems even if the CIO's cannot be easily obtained a priori on both ends.

## *Initial Condition Based on Conjugate Formulation*

While the integration process (9) has fixed boundary points for $\lambda$ from $\pm\infty$ to 0, in practice it can be awkward or numerically tricky to define the derivatives of $k_1, k_1, k_1, \ldots k_N$ w.r.t. $\lambda$ at large enough values of $\lambda$ to mimic $\lambda = \pm\infty$ at the exact CIO for $h$. One can alternatively start this integration process based on the alternative formulation (12) with $\mu = 1/\lambda$ being the integration variable, and start with well-behaved derivatives of $k_1, k_1, k_1, \ldots k_N$ w.r.t. $\mu$ at $\mu = 0$ according to (12). This leads to a robust recipe to start the integration process at the exact CIO of $h$ without need to artificially displace the starting point, which can be switched back to $\lambda$ based integration at any later point. The solution path of the tradeoff curve can then be viewed as having fixed starting point $\mu = 0$, and fixed end point $\lambda = 0$ (See (20) below).

## *Global Optima*

When the system is nonlinear, multiple local optima can exist for a given constraint. It is beyond the scope of this note to analyze the structure of these optima and its implication for the algorithm discussed. It is however observed, as is demonstrated in the examples below, that the deterministic solution path traverses the same constraint (or objective depending on the view of the problem) multiple times on the trade-off landscape, each time producing a local optimum for that given constraint. This provides a collection, whose completeness awaits proof, of local optima from which a quasi-global optimum can be extracted. It is conjectured here that this may be the true global optimum. More advanced methods would be required to determine if this gives the true global optimum.

---

[14] In practice a sufficiently large $\lambda$ is a good approximation. Alternatively one can use the conjugate formulation to start the process. See below.

[15] Actually one can start from a solution with $\lambda$ equal to any value if it can be known a priori. The CIO is special case for $\lambda = -\infty$.

The question of true global optimum may be related to the uniqueness of the tradeoff curve for a given set of competing objectives, which was implicitly assumed. Given the uniqueness of the CIO and the deterministic nature of the integration algorithm, the reality may be favorable to the above conjecture.

## *A Numerical Example*

In the following the problem of "betatron matching" in magnetic optics is used to illustrate the points made earlier. Betatron matching amounts to using a set of focusing magnets, or quadrupoles, to shape charged particle beams to desired shape in accelerators. A total of $N_Q$ quadrupoles can be set to strengths $k_1, k_2, k_3, \dots k_N$, which impact the input beam shape and cause the output beam shape to change. A matching criterion can be defined through the so called 4D mismatch factor, basically measuring how far the output Real beam shape deviates from the Design beam shape. It will be identified with the symbol $F$ in the context of the current notation.

$$F = \frac{1}{4} Tr\left(\Sigma_D^{-1} \cdot M(k_m) \cdot \Sigma_R \cdot M^T(k_m)\right), \quad m = 1, 2, \dots N_Q$$
$$\Sigma^{ij} = \frac{1}{n} \sum_{k=1}^{n} x_k^i \cdot x_k^j \quad i, j = 1,2,3,4 \tag{18}$$

The 4 coordinates for each of the $n$ particles in the beam are its two transverse positions and two transverse angles. $\Sigma_D$ ($\Sigma_R$) represents the Design (Real) 4×4 beam covariance matrix, and $M(k_m)$ is the 4×4 beam transport matrix as function of the $N_Q$ quadrupole strengths. $M(k_m)$ can take on complicated forms for complicated systems. The goal is to minimize $F$, which grows as matching deteriorates, and equals 1, its lower bound, when the beam is fully matched.
The competing objective $H$ is the quadratic sum of deviation of Real quadrupole strengths from Design, with an interest of being minimized, technically speaking, to its lower bound 0:

$$H = \sum_{m=1}^{N_Q} (k_m^R - k_m^D)^2 = \sum_{m=1}^{N_Q} \delta k_m{}^2 \tag{19}$$

The CIO for $H$, or (6), is the trivial solution $\delta k_m = 0$ for all $k_m$. This provides an excellent point to start the program at $\mu = 0$ ($\lambda = -\infty$). The derivatives of $k_m$ w.r.t. $\mu$ at $k_m = 0$ can be conveniently obtained from the conjugate formulation (12), from which the entire integration is launched:

$$\boldsymbol{N_{ij}}\Big|_{\mu=0} = \frac{\partial^2 H(\boldsymbol{k})}{\partial k_i \partial k_j} = 2\delta_{ij} \quad \rightarrow \quad \frac{d\boldsymbol{k_i}}{d\mu}\Big|_{\mu=0, k_m=0} = \frac{1}{2}\frac{\partial F(\boldsymbol{k})}{\partial k_{\boldsymbol{i}}}\Big|_{k_m=0} \tag{20}$$

The CIO for $F$, or the best achievable matching for the system at hand, is unknowable a priori. It is the objective of accelerator control, and can be a source of operational inefficiency and uncertainty if an inadequate numerical algorithm is used for its solution. It will be shown that the program presented here will start from the CIO for $H$, and reach the CIO for $F$ when $\lambda$ becomes 0 following the deterministic path.
Specifically the example is shown of a 6-quadrupole system at 15° betatron phase advance intervals. The physical meaning of this is that these quadrupoles have poor orthogonality characteristics in terms of their ability to collectively control beam shape, thus making a numerical solution more labored and unpredictable. It is very conceivable that straightforward application of a naïve numerical algorithm to (18) and (19), as well as any tradeoff thereof, would lead to compromised solutions. Further computational detail and outcome of this program applied to the 6-quadrupole betatron matching is discussed below[16].

Direct comparison to typical solver and algorithmic detail

As control case the numerical tool *Mathematica* [2] was used to solve for the CIO's for $F$. There are a few options in its repertoire: *Solve*, *NSolve*, and *FindRoot*. For the computer used for this test[17], *Solve* and *NSolve* entered a prolonged state of no response with some appreciable usage of RAM, and were aborted after a couple of hours. *FindRoot* required an initial guess, for which $\delta k_m = 0$ was used for all $k_m$ for lack of better guidance, and was able to produce results after reasonable computing time. *FindRoot* was able to achieve the true CIO for $F$, which is 1 in this 6-quadrupole system, in many cases, and failed to do so in some other cases, including ones where the shortfall is considerable. It is understood that given a better initial guess *FindRoot* would doubtless have found the right CIO. In a sense the program proposed in this note serves exactly that purpose.

In contrast the proposed program found the true CIO in all cases[18]. Instead of actually carrying out the integration (9) and (15), which required software resource not readily available, we mimicked this process by Alternating *FindRoot* on (1) and (3) repeatedly in tiny increments, advancing along the solution path, taking as initial guess what was generated in the previous step. The direction of increment, in $h, f$, or $\lambda$, is dictated by the previous step as well. When *FindRoot* fails to produce an exact solution, it is indication that an inflection or extremum is encountered, and an iterative procedure is invoked to locate the exact position of the inflection or extremum before switching over to the complementary process. This is equivalent to the more elegant approach by (9) and (15), which requires serious software development not available at this point. The program can take considerable time to complete when the number of quadrupoles exceeds 4 or 5 and the betatron phase advance narrows. This is however predominantly, likely at over 99% level, due to the time spent by *FindRoot* in iteration when the exact solution is not achievable at an inflection or extremum, as mentioned above. This situation is exacerbated in the iterative phase. With careful tuning of the process parameters this execution time might be minimized considerably. There is no effort taken on this fine tuning for the current proof-of-principle program[19].

Perhaps the power of the current approach is best demonstrated in some examples with 4 quadrupoles. By degree-of-freedom this is the critically constrained situation where the number of variables ($\delta k_m$) is exactly enough to satisfy the condition of mismatch factor $F$=1 for arbitrary input condition. This however does not guarantee that the solution is real, and in some cases only complex solutions exist [2]. In such cases the current approach unambiguously reaches $\lambda$=0 and terminates, with $F$ stopping at a value greater than 1, which is the best matching achievable with real $\delta k_m$. On the other hand relying on a numerical root-finding or optimization algorithm cannot unequivocally answer the question of whether the best possible solution has been reached, or even what this solution might be. The same applies to cases with less than, and in some rare cases more than, 4 quadrupoles.

---

[16] In many pictures shown here, for technical reason, the objective $F$ takes on the symbol $CS$, and the constraint $H$ the symbol $K$.
[17] Intel i7-3820 CPU @ 3.60 GHz, 64 GB RAM, 64-bit OS.
[18] An elaborate process was carried out to scan the 4D input space of all possible beam distortions.
[19] In addition this application is aimed at building an interpolation database for fast online accelerator control. Thus offline generation of the database is not particularly time sensitive.

It is also notable that the program required almost no tuning over the wide range of input conditions explored. Scanning over input parameters was run in automatic mode, repeating this process from one case to the next without user intervention. It is however conceivable some parameters in this process may warrant tuning when the class of problem is changed.

Another comparison was made with a common practice in such problems where $F$ and $H$ are combined, often quadratically, to form a single objective for optimization. An artificial weighting factor between $F$ and $H$ must be introduced. This can be visualized in Figure 5. Again we used *FindRoot* as the algorithm for solution, and it reproduced the correct optima on the tradeoff curve with hit or miss at a level similar to the test before. Very often the arbitrariness in choosing the weighting makes it difficult to either control or interpret the optimized outcome, and making this approach of less value in general.

A typical solution path and reason for insisting on $\lambda$=0 as termination criterion

Figure 6 gives the detailed structure of a typical solution path in the $\lambda$–$f$ plane obtained by alternating (5) and (9), after starting $\lambda$ off at a large negative value[20]. The only criterion for process termination is $\lambda$=0 ($f$=1). At several points in the process this criterion was all but exactly met ($\lambda$=–2.42; $f$=1.0076 and $\lambda$=–0.74; $f$=1.00013), but the program continued, in the second case taking a wide loop away from the optimum before returing to the absolute optimum at the end. Table 2 lists the underlying solutions for $k_m$ and $\lambda$ around the first turn-around defined by an inflection and an extremum of $f$. These two points are numerically very close and presents a measure of resolution the integration method must meet. This does not pose a problem for the proof-of-principle program currently used.

Besides being rigorous and insisting on reaching the true optimum as opposed to approximate ones, there is a practical advatnge in adhering to rigor. This is seen in Figure 7 where the same path is shown in the $h$–$f$ plane, or the plane of objective ($CS$: mismatch factor) and constraint ($K^2$: magnet strength) in the physical picture. Had we stopped at the approximate optimum ($\lambda$=–0.74; $f$=$CS$=1.00013) we would have obtained a solution for $h$ ($K^2$) of 0.008. By inisiting on following the program to the true optimum, with insignificant gain in either $\lambda$ or $f$, $h$ ($K^2$) is reduced to 0.003, a major gain.

The solution for each of the 6 quadrupole strengths along the solution path is shown in Figure 8.

Point of diminished return

Correspondence between extrema in the $\lambda$–$h$ tradeoff curve and inflections, or points of diminished/enhanced return, in the $h$-$f$ tradeoff curve discussed earlier in (11) is emphatically underscored by some examples seen in the numerical solutions. One example is given in Figure 2(B). After making substantial gain (2.5 to 1.05) in reducing the mismatch factor ($F$) at little cost (~0.004) to quadrupole strength ($H$) up to the inflection point, it would take much more (>0.008) quadrupole strength to gain only another 0.05 in mismatch factor. This point is again made in Figure 13.

Global optimal tradeoff

Every point on the solution curve represents a locally optimal tradeoff between $f$ and $h$. Figure 7 however reveals multiple such locally optimal $f$ at many given $h$, and vice versa. In the specific example given where both $f$ and $h$ are to be minimized, the lowermost and leftmost braches of the curve give the globally optimal tradeoff. This works for both $f$ and $h$ unambiguously, because the slope of this curve is everywhere negative by definition. One can thus take the two parts of the curve intersecting near $f(CS)$=1.02, disconnect and join the two lowest/leftmost sections at the intersection and define this spliced curve as the global optimal tradeoff. Everything above and to the right is sub-optimal, as shown in Figure 9. This works simultaneously and unambiguously for both "minimal $f$ for given $h$" and "minimal $h$ for given $f$", because the slope is everywhere negative[21]. At the level of actual solution for the $k_m$'s, this amounts to splicing and rejoining each curve in Figure 8 at $f(CS)$=1.02, as also demonstrated in Figure 9.

[20] In practice the scheme mimicking (5) and (9) using *FindRoot* was used as discussed.

The process of splicing the local tradeoff curve to obtain the global optimal tradeoff is akin to that of extracting the "Pareto Front" [3] in multi-objective genetic optimization practice. As discussed earlier, the distinction between objective and constraint is artificial and we are really dealing with a multi-objective problem. The Global Optimal Tradeoff curve in Figure 9 is Pareto-dominant over the rest of the local curve, and thus constitutes the global optimum with the following features:

- Of all points having the same $f$ value, the one on the Global Optimal Tradeoff dominates all others in terms of $h$ value, and vice versa.
- The Global Optimal Tradeoff curve is monotonic, despite the twists and turns in the local tradeoff curve.

Both features depend strongly on the fact that λ has to be negative everywhere.

Figures 10 and 11 show the global behavior of the solution curve in terms of λ, $f$, $h$, and all $k_m$'s. This gives some insight into the singular behavior and nature of "turn-around" points discussed earlier.

## *A Second Numerical Example*

Instead of fixing a mismatch in the input beam as in the previous example, this algorithm is also applied to a variation of the problem where the input beam is matched, but the optical transport itself is in error and results in final mismatched beam. The formulation is almost the same except that the error introduced is no longer minor, but caused an initial mismatch factor, eqn.(18), of over 7000! This is a good test of the robustness of the program. Furthermore we have the option now of not only minimizing the RMS $\delta k_m$ (eqn. (19)), but also RMS $k_m$ itself since the baseline optics design is no longer valid anyway.

$$H = \sum_{m=1}^{N_Q} {k_m}^2 \tag{21}$$

Some of the results are shown in Figures 12 and 13. The observation is that the algorithm is robust enough that, without further tweaking of parameters, it brought the mismatch factor Φ from over 7000 to 1. The interplay between minimum in the λ–$f$ curve and the point of diminished return as indicated by contrast between quadrupole strength and mismatch is also evident in the last example in Figure 13.

More cases were tested by placing the optical transport error at different locations. This led to the results shown in Figure 14. It is interesting to note that in some cases choosing constraints of minimal RMS $\delta k_m$ or minimal RMS $k_m$ led to different final solutions, although the latter should be formally constraint independent. It is also worth noting that in each case the final solution is consistent with the constraint chosen. This strongly suggests more complex topological structure of the solution space in at least some of the cases, as further speculated in Figure 16.

---

[21] This will still be true with problems of maximization etc. as $\lambda$ will change sign accordingly.

It is interesting to note that the notion of Globally Optimal Tradeoff extracted via a Pareto domination procedure still holds in each scheme's respective context. For example if solution A is derived via the constraint $h$ on minimal RMS $\delta k_m$, then even though there exists another solution B with exactly the same optimal objective $f$ (CS=1) derived via minimal RMS $k_m$, B is still Pareto-dominated by A because B has larger RMS $\delta k_m$, and vice versa. The algorithm settles into different Globally Optimal Tradeoff curves consistent with the constraint chosen. The Pareto front contains much more information than end points alone.

## *Alternative Singularity-Free Formulation for Numerical Integration*

Although the procedure outlined in (9) through (15) is realizable via a Runge-Kutta type implementation, the need to circumvent singularities may present numerical issues. As mentioned, numerical results presented so far were obtained by actually scanning the solutions to (1), (2), or (3) in micro steps, which is equivalent to a Runge-Kutta type integration but does not guarantee a smooth execution of the latter. In addition this approach is only for proof-of-principle, lacking features essential to large scale, efficient application, such as adaptability in step size. Furthermore this stepwise solution approach takes advantage of a priori knowledge of a one-dimensional "direction" of integration, either in $f$, $h$. or $\lambda$. We steer the solution in the next step to be in the direction of increasing or decreasing single parameter. This would be difficult to realize at the level of, for example, (9) or (15), where one cannot know a priori what the directions for the many $k_m$'s should be in the next step. The above motivates the formulation of the integration in a more robust manner free of singularity. This is done by going to the independent variable $|dk|$, the length along the tradeoff curve in the $\boldsymbol{k}$–space.

It is easily seen by squaring (9) ,

$$\begin{aligned} |dk|^2 &= d\lambda^2 \frac{Adj(\boldsymbol{M}) \cdot \boldsymbol{R} \cdot \boldsymbol{R}^{\boldsymbol{T}} \cdot Adj(\boldsymbol{M})^{\boldsymbol{T}}}{Det(\boldsymbol{M})^2} \\ &= \frac{d\lambda^2}{Det(\boldsymbol{M})^2} \cdot |Adj(\boldsymbol{M}) \cdot \boldsymbol{R}|^2 \end{aligned} \tag{22}$$

that the direction cosines along the tradeoff curve are

$$\frac{d\boldsymbol{k}}{dk} = \pm \widehat{\boldsymbol{P}}, \qquad \boldsymbol{P} = Adj(\boldsymbol{M}) \cdot \boldsymbol{R} \tag{23}$$

where we shorthanded $|dk|$ as $dk$, bold faced letters $\boldsymbol{k}$ and $\boldsymbol{P}$ denote vectors, and carets denote unit vector. Sign ambiguity introduced when squaring in (22) is necessary and should be resolved by choosing the direction such that $\boldsymbol{k}$ advances in the desired direction[22]. Some deliberation shows that $\boldsymbol{P}$ can never be singular, or even a zero vector barring some pathological $\boldsymbol{M}$ whose example is difficult to come by[23]. Likewise one can show that

$$\begin{aligned} \frac{d\lambda}{dk} &= \pm \frac{Det(\boldsymbol{M})}{|\boldsymbol{P}|} \\ \frac{df}{dk} &= \pm \frac{\lambda \cdot (\boldsymbol{R}^{\boldsymbol{T}} \cdot Adj(\boldsymbol{M}) \cdot \boldsymbol{R})}{|\boldsymbol{P}|} = \pm \lambda \cdot \boldsymbol{R}^{\boldsymbol{T}} \cdot \widehat{\boldsymbol{P}} \\ \frac{dh}{dk} &= \pm \frac{(\boldsymbol{R}^{\boldsymbol{T}} \cdot Adj(\boldsymbol{M}) \cdot \boldsymbol{R})}{|\boldsymbol{P}|} = \pm \boldsymbol{R}^{\boldsymbol{T}} \cdot \widehat{\boldsymbol{P}} \end{aligned} \tag{24}$$

---

[22] More precisely, depending on whether one is going from maximum (minimum) to another extremum of the same or opposite flavor, this sign is uniquely fixed.

[23] Basically inverse of a normal matrix $\boldsymbol{M}$ cannot be rank deficient since by definition its own inverse is just $\boldsymbol{M}$.

The signs chosen in (24) should be the same as that in (23). These can be shown to be consistent with (9) and (11). Again all these expressions are free of singularity. The first equation can be 0 with $\boldsymbol{M}$ being rank deficient, and the second and third equations can be 0 with $Adj(\boldsymbol{M})$ (or $\boldsymbol{M}^{-1}$) satisfying the saddle point condition of having both positive and negative eigenvalues. A zero in any of the above equations marks a "turn-around" point in the previous formulation without introducing singularity in the equation to be integrated, (23), which never stalls by virtue of $\widehat{\boldsymbol{P}}$ becoming 0 either. The process stops unambiguously when λ=0.

Since $\boldsymbol{M}$ depends on $\lambda$, the Runge-Kutta integration must evaluate both $\lambda$ and $\boldsymbol{k}$ at each step, but we are free from the need to know which way *f*, *h*. or $\lambda$ needs to go in the next step. Everything is dictated by well-behaved local derivatives.

Concern to avoid initial infinite $\lambda$ leads to the conjugate expression with $\mu = 1/\lambda$ as before:

$$\frac{d\boldsymbol{k}}{dk} = \pm\widehat{\boldsymbol{Q}}, \quad \boldsymbol{Q} = Adj(\boldsymbol{N}) \cdot \boldsymbol{S}$$
$$\frac{d\mu}{dk} = \pm\frac{Det(\boldsymbol{N})}{|\boldsymbol{Q}|}$$
$$\frac{df}{dk} = \pm\frac{(\boldsymbol{S}^{\boldsymbol{T}} \cdot Adj(\boldsymbol{N}) \cdot \boldsymbol{S})}{|\boldsymbol{Q}|} = \pm\boldsymbol{S}^{\boldsymbol{T}} \cdot \widehat{\boldsymbol{Q}} \tag{25}$$
$$\frac{dh}{dk} = \pm\frac{\mu \cdot (\boldsymbol{S}^{\boldsymbol{T}} \cdot Adj(\boldsymbol{N}) \cdot \boldsymbol{S})}{|\boldsymbol{Q}|} = \pm\mu \cdot \boldsymbol{S}^{\boldsymbol{T}} \cdot \widehat{\boldsymbol{Q}}$$

It is important to avoid integration variables with disparate magnitudes in a Runge-Kutta procedure to ensure numerical stability. Therefore the best strategy is to start with (25) at μ=0, integrate until μ = λ = –1[24], then switch to (24) and continue integrating until λ=0. One can even further scale the system such that all $k_m$'s vary within a range of order unity. This would improve the overall numerical stability.

It is interesting to observe that (23) and (25) indicate that, at the end of the optimal "constraint" (*H*), it is the form of the "objective" (*F*) that is driving the course of the tradeoff curve, and vice versa. This is necessary or the process would stall right at the start. Thus at every optimal point (5) or (6) there are an infinite number of tradeoff curves emanating from it, driven by different competing objectives. This picture is intuitively sensible. It is still an outstanding question, lacking deeper understanding of the topological structure of the solution space, how a formally constraint-independent optimum (5) or (6) can only be accessed from one constraint but not the other, as suggested in Figure 14 and speculated in Figures 16 and 17. This may be hardly surprising as $\nabla F = 0$ is only a local statement, saying nothing about the global topology of the solution space. Neither is the relationship obvious between two Tradeoff Curves linking the same objective on one end to different competing objectives on the other end, although such situation must exist by virtue of the recipe of integration. All these cannot be answered without a deeper understanding beyond the scope of the current work.

At the end of this process when the entire tradeoff curve is obtained, it is still necessary to look at the *f*–*h* diagram and extract the Global Optimal Tradeoff curve via a Pareto-domination type procedure, which in turn dictates the splicing of the $k_m$ solution curves into global optimal solutions.

[24] Strictly speaking this only applies to the special case of going from maximum (minimum) of one objective to the maximum (minimum) of another objective with no intervening extrema. In more complicated landscape, as seen in the following examples, one can encounter the need to switch more frequently in cases where the tradeoff characteristic between the two objectives changes sign, namely, $\lambda$ or $\mu$ itself changes sign half way. In such cases switching may be warranted to avoid singular situation, although again this is not something fundamental.

## *Numerical Examples Based on Integrating New Formulation*

Examples given below are results of actually integrating the differential expressions (23), (24) and (25)[25] instead of performing local solutions of the Lagrange condition in miniature steps as in some earlier examples. Two examples are given.

Example 1

The two objective functions are

$$F(x,y) = (x - x_0)^2 + (y - y_0)^2$$
$$H(x,y) = \frac{x^2}{a^2} + \frac{y^2}{b^2}$$

It is trivial to derive the parametric expressions for the tradeoff curve analytically

$$\big(x(\mu), y(\mu)\big) = \left(\frac{a^2 x_0 \mu}{-1 + a^2 \mu}, \frac{b^2 y_0 \mu}{-1 + b^2 \mu}\right), \qquad 0 > \mu > -1$$
$$\big(x(\lambda), y(\lambda)\big) = \left(\frac{a^2 x_0}{a^2 - \lambda}, \frac{b^2 y_0}{b^2 - \lambda}\right), \qquad -1 < \lambda < 0$$

Figure 18 demonstrates how integrating (23), (24) and (25) exactly reproduced these curves between the two optima.

Example 2

The two objective functions are periodic in x & y:

$$F(x,y) = \mathrm{Sin}[px] + \mathrm{Sin}[qy]$$

$$H(x,y) = \mathrm{Sin}[rx] + \mathrm{Sin}[ty]$$

The formula for the tradeoff curve can be written in closed form but only solved numerically:

$$\frac{p \cdot t}{q \cdot r} = \frac{\mathrm{Cos}[rx] \cdot \mathrm{Cos}[qy]}{\mathrm{Cos}[px] \cdot \mathrm{Cos}[ty]}$$

Figure 19 shows the contours of *F* and *H*, as well as the numerically computed tradeoff curves. Due to the periodic nature of the functions, and the fact that now we have both maxima and minima, the contours are much more convoluted, including lines containing saddle points in both *F* and *H*. Many links of the tradeoff curves have been independently computed by integrating (23), (24) and (25) with exact agreement. Due to the presence of both maxima and minima now, when the integration approaches a saddle point, namely, when $\lambda$ or $\mu$ crosses 0, additional switching between $\lambda$ and $\mu$ needs to take place to keep the integration going, but this is more a cosmetic procedure than anything substantive.

[25] Using NDSolve in Mathematica.

## *Summary*

In this note we discussed the following

- A deterministic program capable of mapping out global optimal tradeoff between 2 competing objectives is presented. There is no need for initial guess or random variable search. It only requires numerical integration algorithm meeting resolution demands of the particular system.
- The key to this program is the integration on the $\lambda$ parameter, to be alternated with either one of the objectives. This $\lambda$-based integration is essential to the process because
  - It relies on system independent universal end points $\lambda$=0 and $\mu$=0.
  - It is needed for negotiating around singular points that integration based on neither of the competing objectives can.
  - In case root-finding algorithms are used to mimic the integration process, it presents a simpler numerical problem.
- An artificial constraint scheme is devised to map out the optimal tradeoff even for systems where the constraint independent optimum (CIO) cannot be easily obtained at both ends of the curve. This is done by imposing an artificial & easy constraint at one end and integrating to the other end to obtain the CIO, followed by integrating backwards along the true (difficult) constraint to map out the entire tradeoff curve. This concept is illustrated in Figure 15.
- Exploration by the locally optimal curve in the tradeoff space provides a global view of multiple optima under a given constraint. This leads to the "global" optimal curve by slicing together sections of the local optimal curve and discarding inferior branches thereof. It is conjectured that this process may lead to the true globally optimal tradeoff curve.
- An alternative formulation more suitable for Runge-Kutta type implementation that is free of singularity is developed. Using both conjugate forms of this formulation and proper scaling of the system should ensure good numerical stability of the process. Some numerical examples were explored.
- An example of magnetic optical matching was used to illustrate the above points.

## *Further Considerations*

- Extension to multiple objective/constraint, or really just multiple competing objectives totaling more than 2, is an interesting problem. In the context of Lagrange multiplier this problem is well studied. In the current scheme such extension has not been attempted but no fundamental obstacle seems obvious.
- Since this program maps out the objective-constraint trade-off over the entire range of the constraint, to what extent this result already explores the question of inequality constraints is an interesting question requiring further understanding.
- The recipe presented may have limited applicability for linear systems because the Hessian of (9) or (12) may become rank deficient in case $F$ or $H$ is linear in the $k_m$'s.
- It is beyond the scope of the current note to investigate the rich global algebraic structure inherent in the formulations (9) and (12).
  - Figure 16 attempts some speculation on the topological possibilities and their implication on this algorithm. To probe into the next level of this problem such global structure must be understood.

- Saddle points of $F$ or $H$ have been given no distinction from extrema by the formulation in this note in the interest of simplicity[26]. In some applications extra measures may be needed to make this distinction.
- Figure 17 reveals one global structure based on a realistic betatron matching example by extending the tradeoff curve beyond the first solution to (5): $\nabla F = 0;\ \lambda = 0$, into further solutions. These correspond to cases where increased quadrupole strength ($H$) leads to alternatingly increased and decreased mismatch factor ($F$). There is strong suggestion that the same tradeoff curve may visit <u>all</u> solutions to (5), saddle points included. This however may be due to the fact that $H$ in our example admits only <u>one</u> extremum and no saddle points, thus the tradeoff curve emanating from it must visit all CIO's of $F$. For more complicated functional form of $H$ this is not obvious.

* Actual numerical integration of (9) and (12) can be done for example a la Runge-Kutta.
  - The tangency conditions (1), (2) and (3) should be imposed at each step for consistency.
  - In this implementation the switchover between (1), (2) and (3) can be triggered by monitoring the relative rate of change per step in (9), (12) and (15) to ensure always using the most smooth derivative, instead of switching only when an extremum or inflection is encountered.
  - A more well behaved integration system (23), (24) and (25) is given that may further enhance the stability of the Runge-Kutta procedure.

---

[26] One needs to be more specific about the meaning of saddle points. In the constraint-independent context of (5) or (6), saddle point is a multi-dimensional property concerning all possible orientations of the $\boldsymbol{k}$ vector, which possibility is of interest here. On the other hand if one stays on the tradeoff curve while crossing this saddle point, only one specific orientation of $\boldsymbol{k}$ is being sampled, and in this strict context the system goes through a <u>one-dimensional</u> maximum, minimum, or inflection, as can be seen in Figure 15.

| $\left.\frac{dA}{dB}\right\vert$ | $A{=}f$ | $A{=}h$ | $A{=}\lambda$ |
|---|---|---|---|
| $B{=}f$ | 1 | $\frac{1}{\lambda}$ | $\frac{1}{\lambda}\cdot\frac{1}{\boldsymbol{R}^{\boldsymbol{T}}\cdot\boldsymbol{M}^{-1}\cdot\boldsymbol{R}}$ |
| $B{=}h$ | $\lambda$ | 1 | $\frac{1}{\boldsymbol{R}^{\boldsymbol{T}}\cdot\boldsymbol{M}^{-1}\cdot\boldsymbol{R}}$ |
| $B{=}\lambda$ | $\lambda\cdot\boldsymbol{R}^{\boldsymbol{T}}\cdot\boldsymbol{M}^{-1}\cdot\boldsymbol{R}$ | $\boldsymbol{R}^{\boldsymbol{T}}\cdot\boldsymbol{M}^{-1}\cdot\boldsymbol{R}$ | 1 |

Table 1. Derivative relations between $h$, $f$, and $\lambda$ along the 1D solution curve. Note for example $d\lambda/df$ can also be expressed alternatively as $1/(\boldsymbol{S}^{\boldsymbol{T}}\cdot\boldsymbol{M}^{-1}\cdot\boldsymbol{R})$.

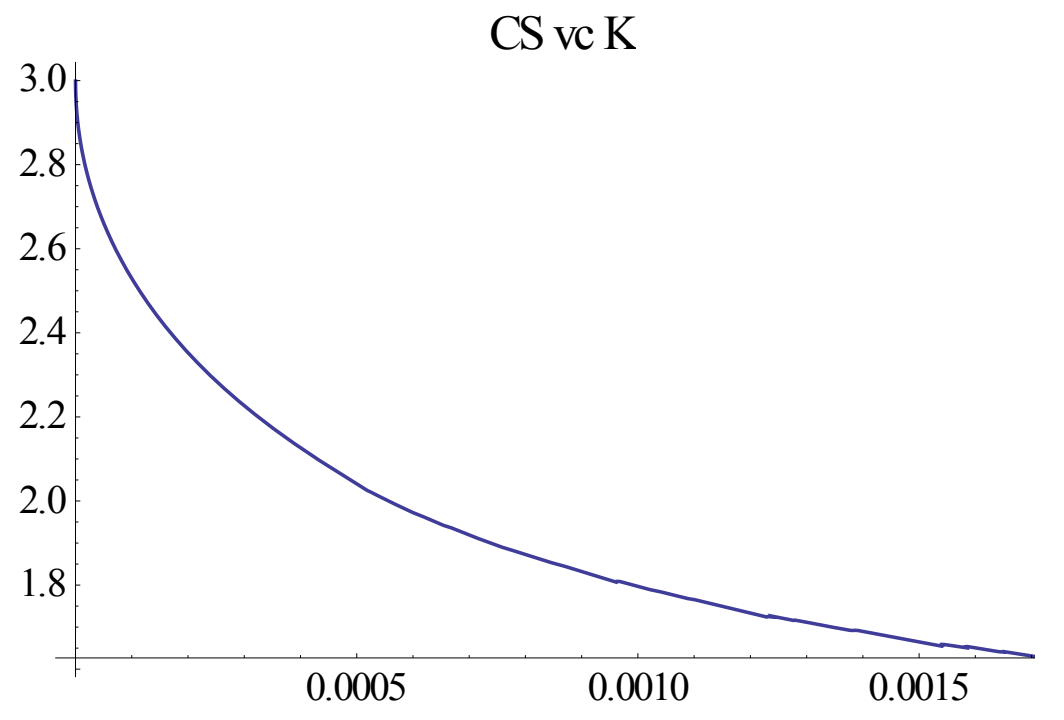

Figure 1. Typical trade-off curve between competing objectives, named *CS* and *K*.

**(A)**

| $\lambda$ vs $h$ | $f$ vs $h$ | |
|---|---|---|
| 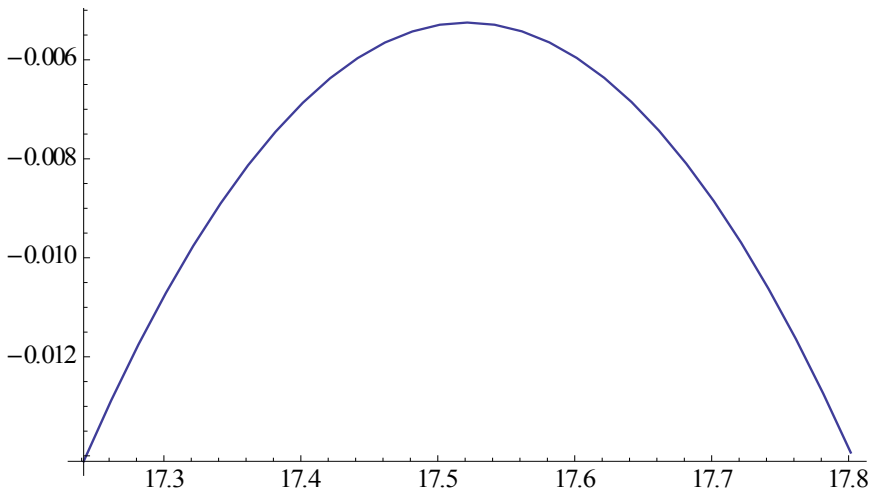  | 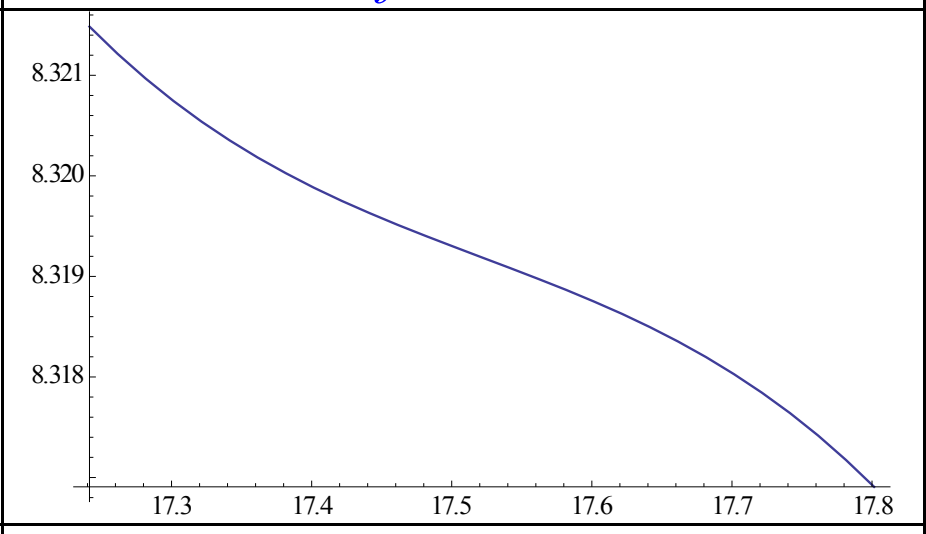  | Matching regime ($\lambda$<0)<br>Point of enhanced return ($\lambda$'=0, $\lambda$''<0) |
| 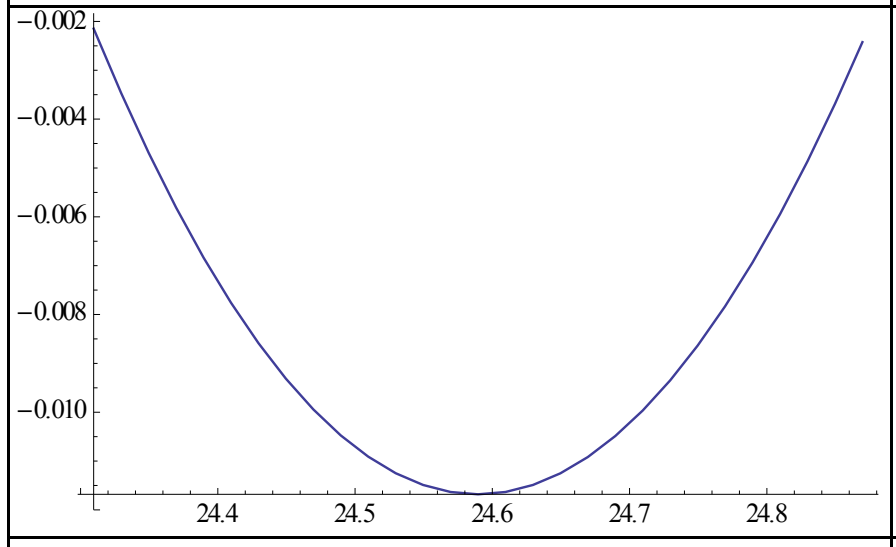  | 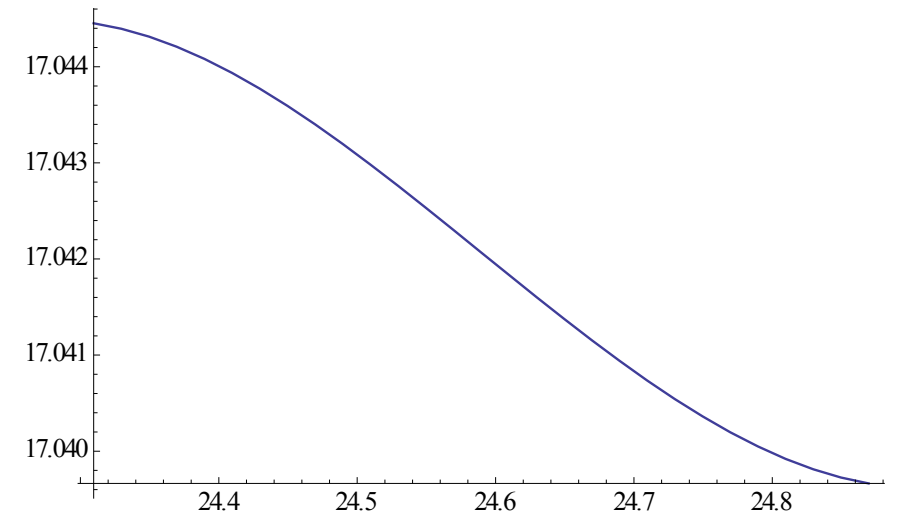  | Matching regime ($\lambda$<0)<br>Point of diminished return ($\lambda$'=0, $\lambda$''>0) |
| 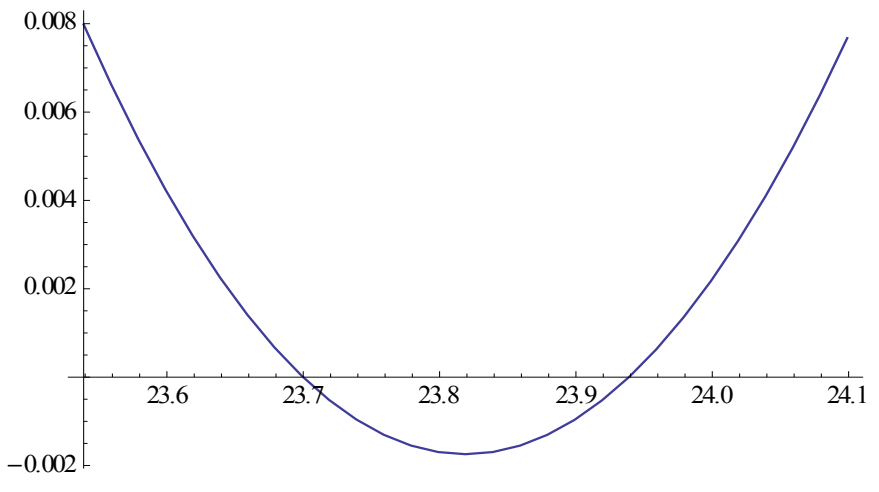  | 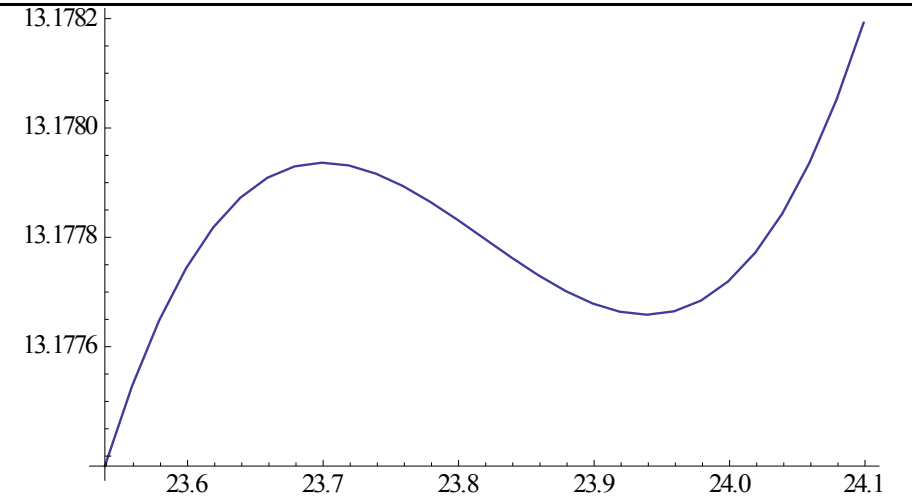  | Matching regime ($\lambda$<0) sandwiched between Anti-matching regimes ($\lambda$>0) |

**(B)**

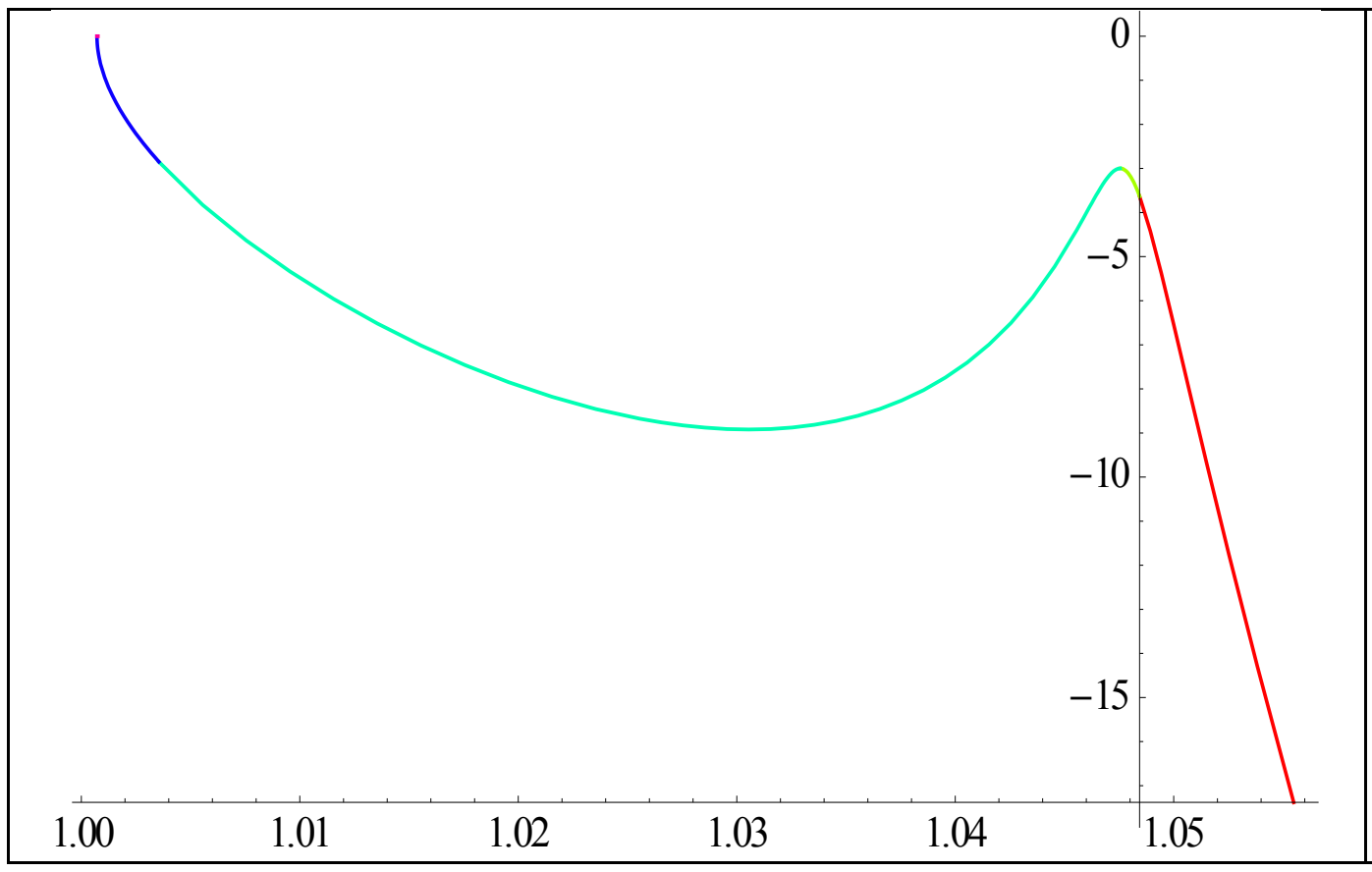

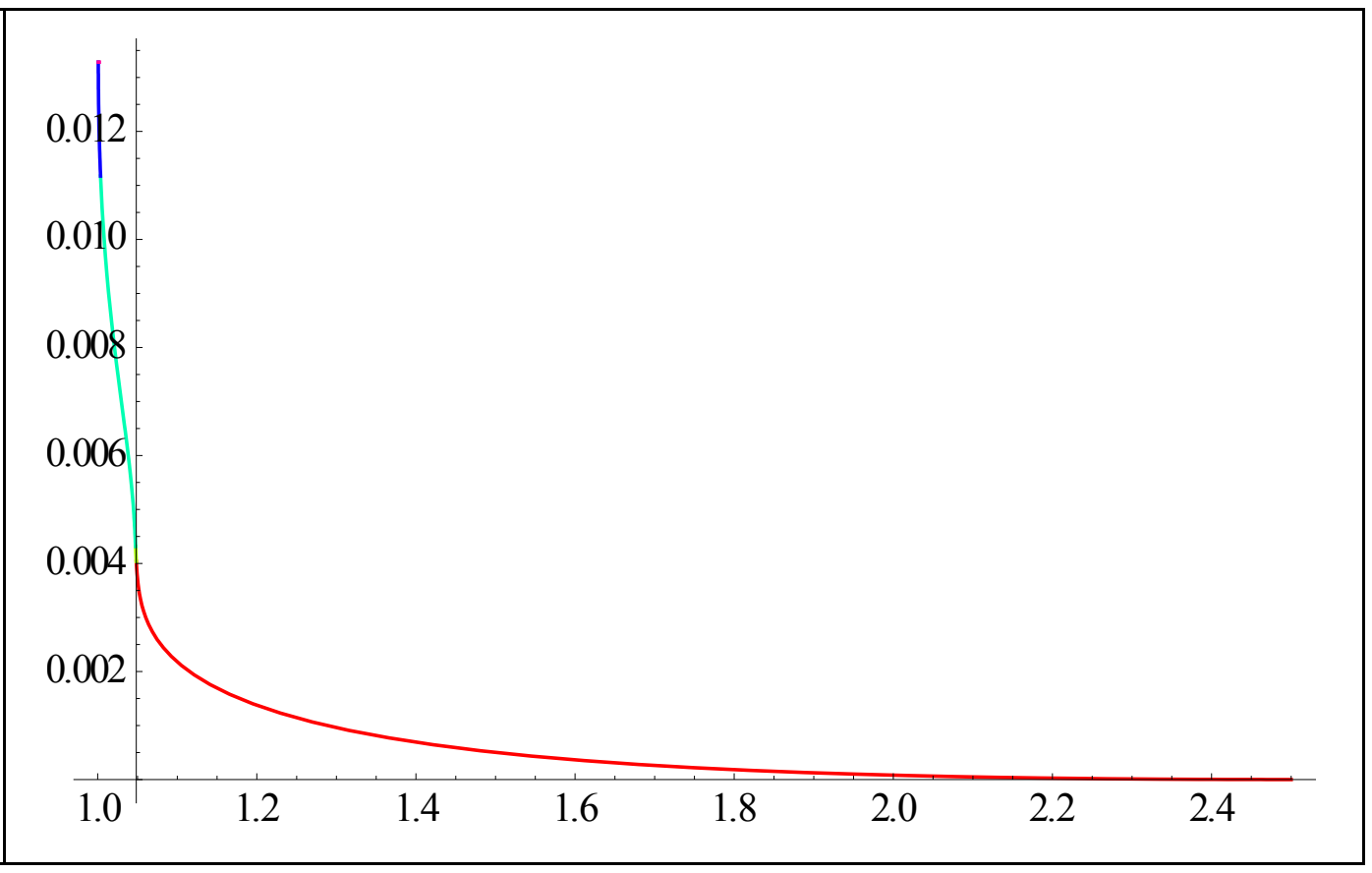

Figure 2. Table **(A)**: Instances of inflection points signified by $Det(\boldsymbol{M}) = 0$ during different stages of constrained optimization taken from a realistic system. Inflection in the $f$ vs $h$ curve can be viewed as points of diminished or enhanced returns beyond which effectiveness of increasing quad strength decelerates or accelerates. Evolution of $\lambda, f$ and $h$ along the solution curve offers a complete picture of the betatron matching process.

Table **(B)**: Inflection point marking pronounced diminished return in realistic betatron matching solution.

Left: $\lambda$ vs mismatch factor ($F$) plot of the end segment of the tradeoff curve right before the final $\lambda = 0$ solution. An inflection in $F$ vs $H$ is clearly indicated by the maximum in this plot near the junction between red and cyan.

Right: Quadrupole strength ($H$) vs mismatch factor ($F$) plot for the entire tradeoff curve from $\lambda = -\infty$ to $\lambda = 0$. Onset of diminished return (gain in $F$ vs expense in $H$) is obvious around the junction between red and cyan.

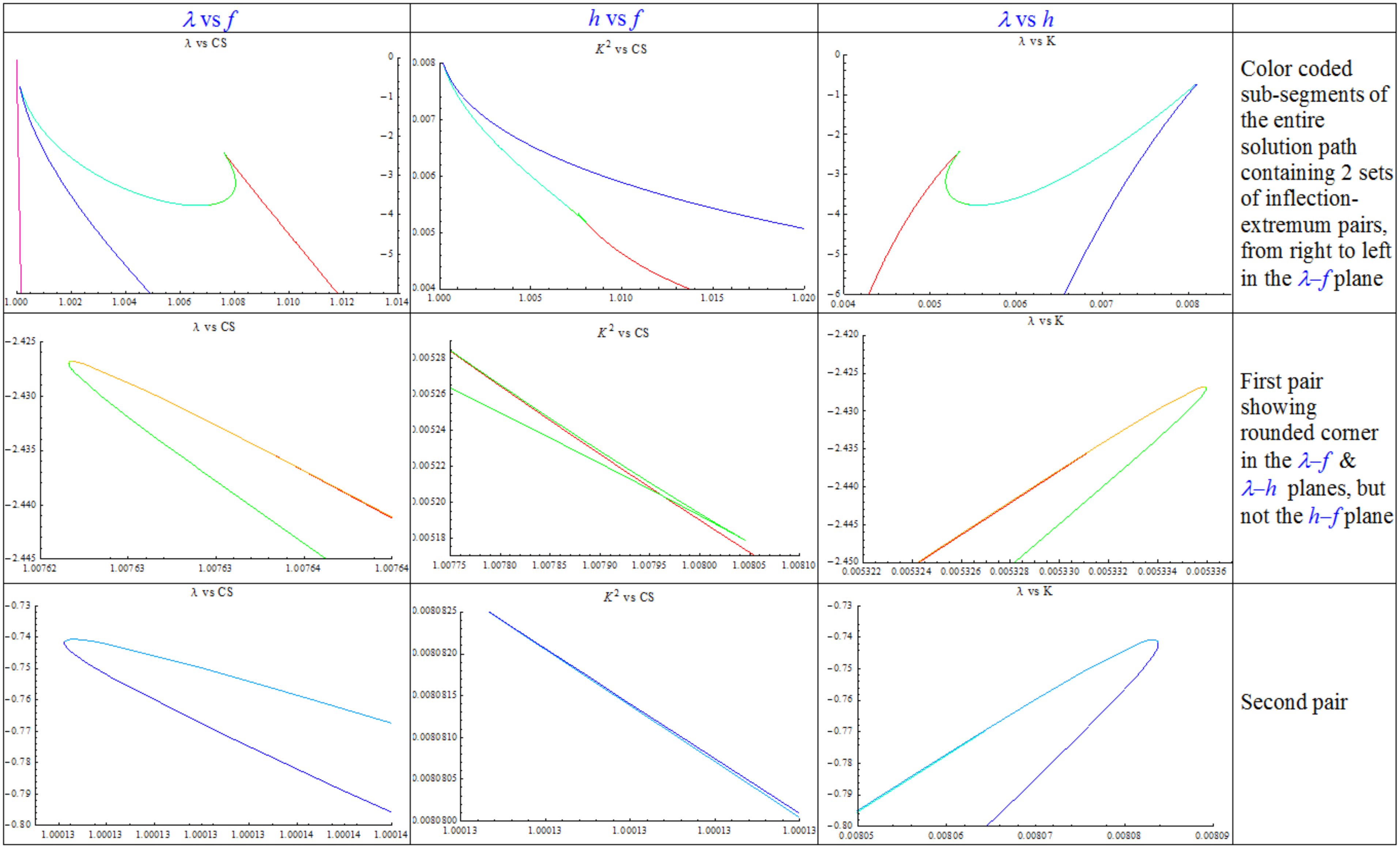

Figure 3. Demonstration of behavior near singular points in either the $\lambda$–$f$, $\lambda$–$h$ or $h$–$f$ plane, taken from a realistic 6-quad betatron matching system. In the former case alternating (1) & (3), or (2) & (3) allowed smooth negotiation around the corner onwards, which is not possible in the latter with only sharp corners. This is because the singularity occurrences of (3), $Det(\boldsymbol{M}) = 0$, are not coincident with those of (1) and (2), $\boldsymbol{R}^{\boldsymbol{T}} \cdot \boldsymbol{M}^{-1} \cdot \boldsymbol{R} = 0$. This is also evident from the fact that the slope of the curve in the $h$–$f$ plane, or $\lambda$, can never change sign within the solution range.

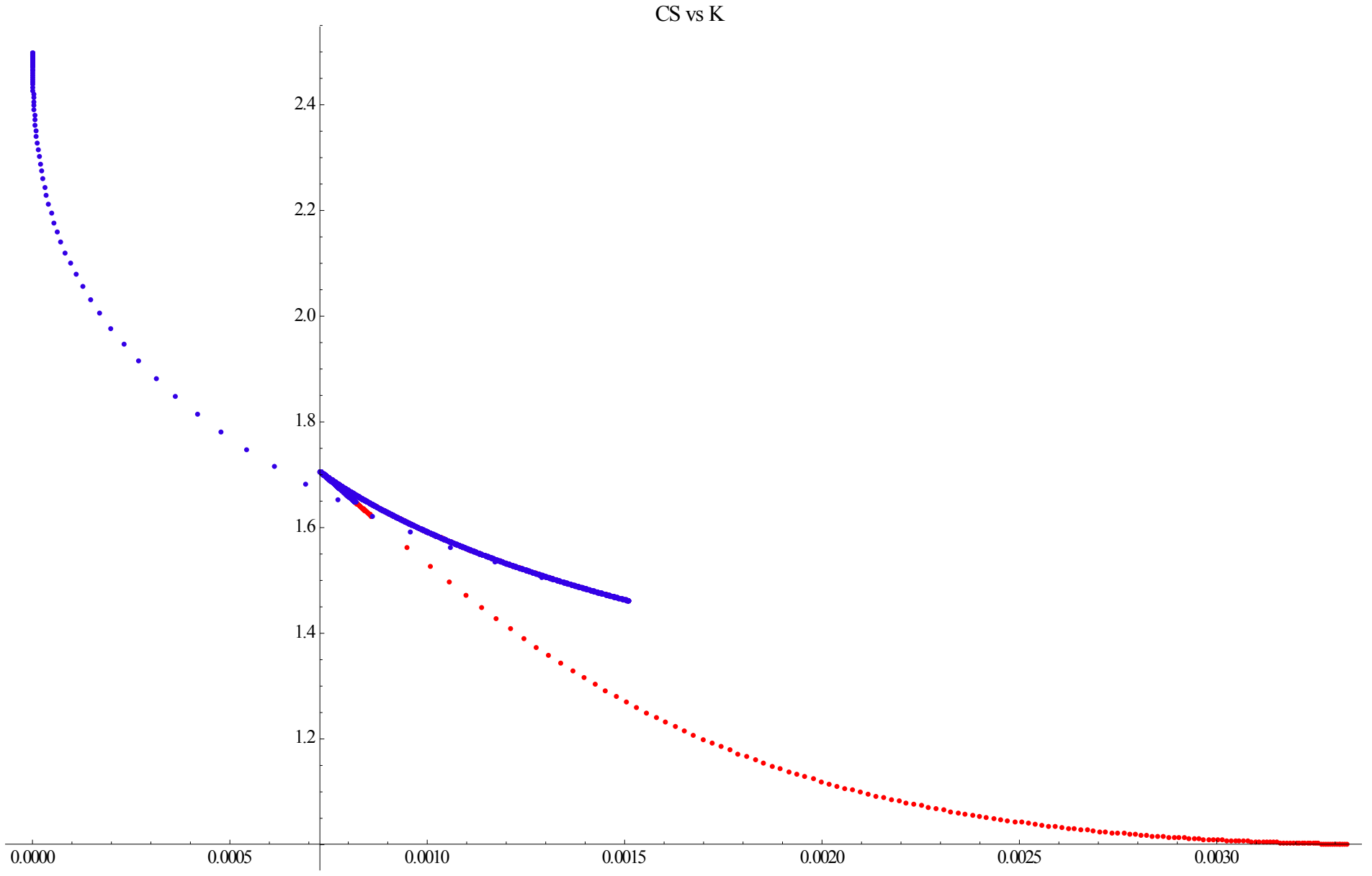

Figure 4. Trade-off curve between objectives $CS$ and $K$, mapped out by integration from the CIO for $K$ (blue) and the CIO for $CS$ (red). The two curves coincide when each is taken to the other end.

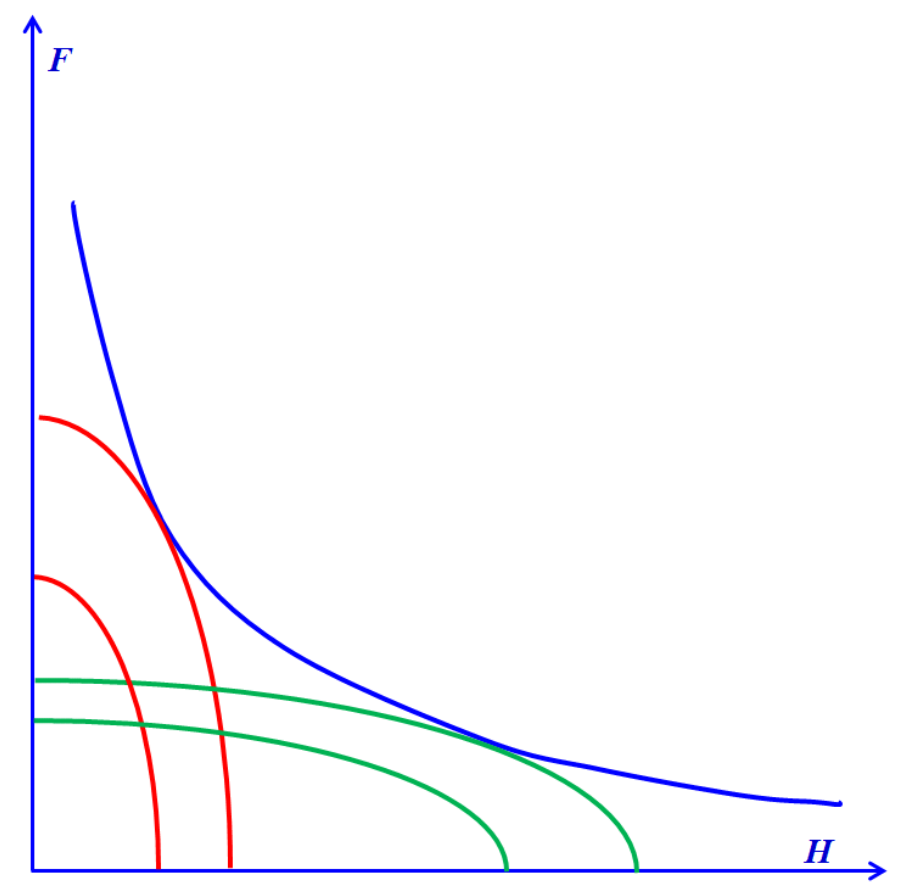

Figure 5. Concept of combining $F$ and $H$ into single objective for optimization. Depending on emphasis on $F$ (green) or $H$ (red), the optimal contour intersects the true trade-off curve (blue) at different locations.

| | Inflection in $F$ | Extremum in $F$ |
|---|---|---|
| $k_1$ | 0.02098529223 | 0.02103888316 |
| $k_2$ | 0.01880655198 | 0.01877436785 |
| $k_3$ | -0.02861344076 | -0.02860804345 |
| $k_4$ | -0.03294511788 | -0.03294056940 |
| $k_5$ | 0.04601841885 | 0.04601961094 |
| $k_6$ | -0.02280337960 | -0.02279314454 |
| $\lambda$ | -2.426805307 | -2.426988991 |
| $\boldsymbol{R}^T \cdot \boldsymbol{M}^{-1} \cdot \boldsymbol{R}$ | $\infty$ | 0 |

Table 2. Solutions for $k_m$ and $\lambda$ around the first turn-around defined by an inflection and an extremum of $f$. $\boldsymbol{R}^T \cdot \boldsymbol{M}^{-1} \cdot \boldsymbol{R}$ went from $\infty$ to 0 within this short stretch.

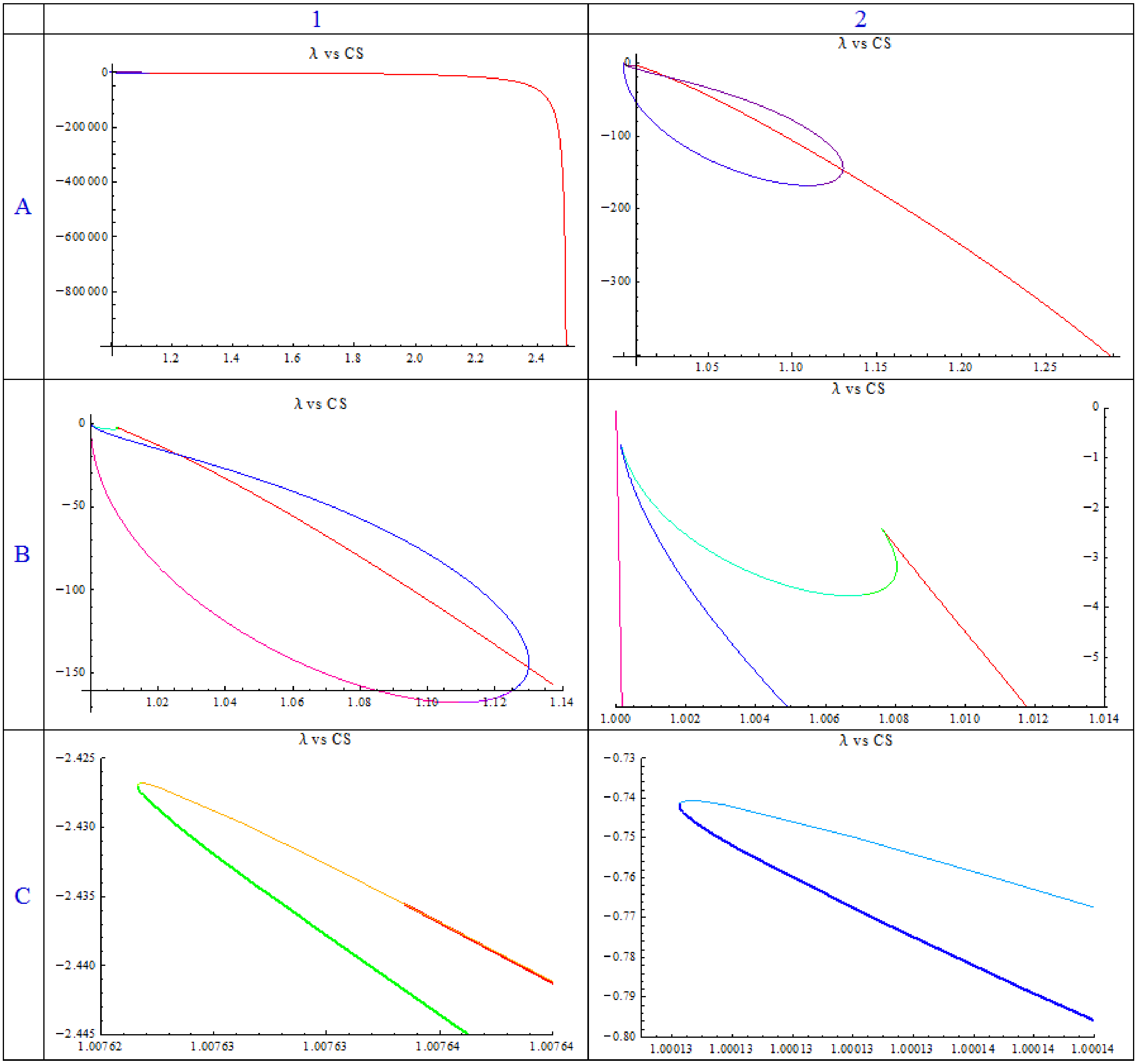

Figure 6. Entire solution path in the $\lambda$–$f$ plane, zoomed to different levels at different locations. Different colors represent different segments of integration.

A1: Global solution path starting at ($f$=2.5, $\lambda \to -\infty$), ending at ($f$=1, $\lambda$=0)
A2: Zoomed in for detail toward the end
B1: Further zoom
B2: Near the point $\lambda$=0 and $f$=CIO=1: The path first follows (5), approaches CIO near 1.0076, encounters an inflection in $f$, switches to (9), encounters an extremum in $f$, switches to (5) again, then again to (9) followed by (5) reaching another inflection at $f$=1.00013 (!), then (9) and (5), executing a loop all the way down to $f$=1.14 before turning around and reaching the true CIO $f$=1, shown by the red line on the extreme left.
C1: Zoom in around $f$=1.0076, the turn-around on the right in B2
C2: Zoom in around $f$=1.00013, the turn-around on the left in B2

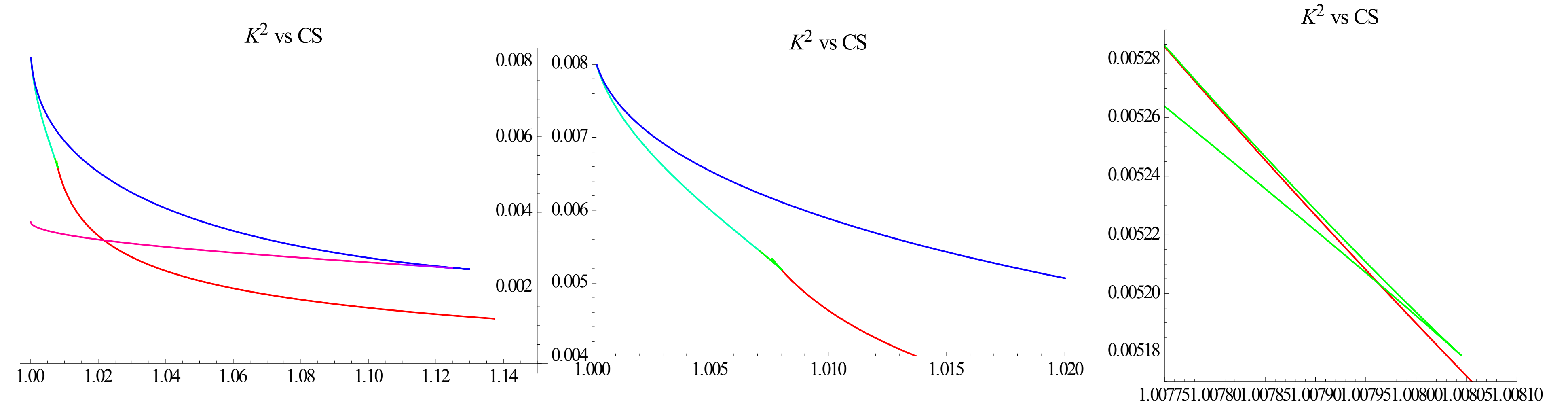

Figure 7. Solution path in the *h–f* plane, zoomed to different levels at different locations. Different colors represent different segments of integration.
Left: Section of path roughly corresponding to B1 in Figure 6 containing all turn-around areas.
Middle: Zoomed in.
Right: Further zoom-in reveals a send bow-tie structure.

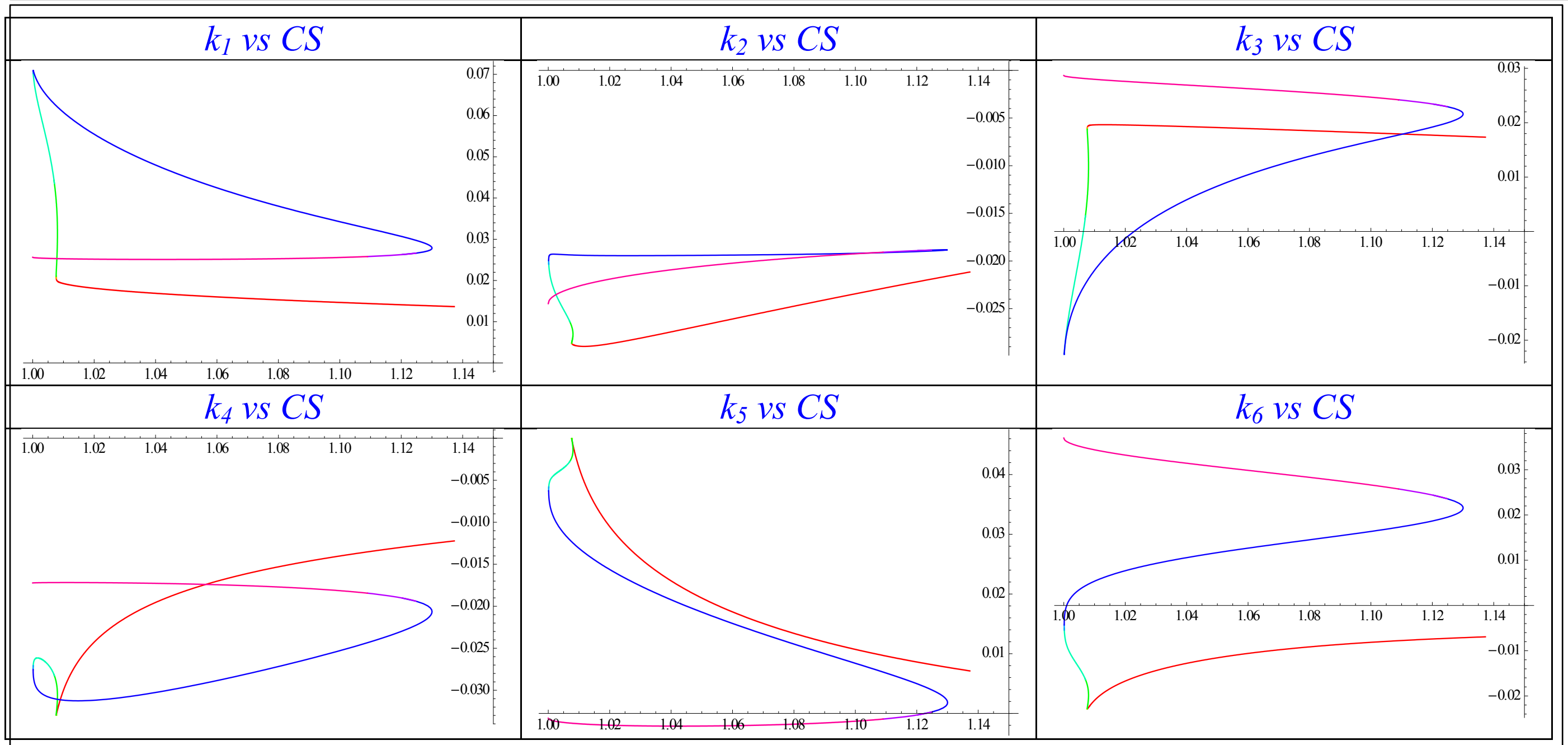

Figure 8. Solution path followed by each one of the 6 quadrupoles, zoomed to section roughly corresponding to B1 in Figure 6 containing all turn-around areas.

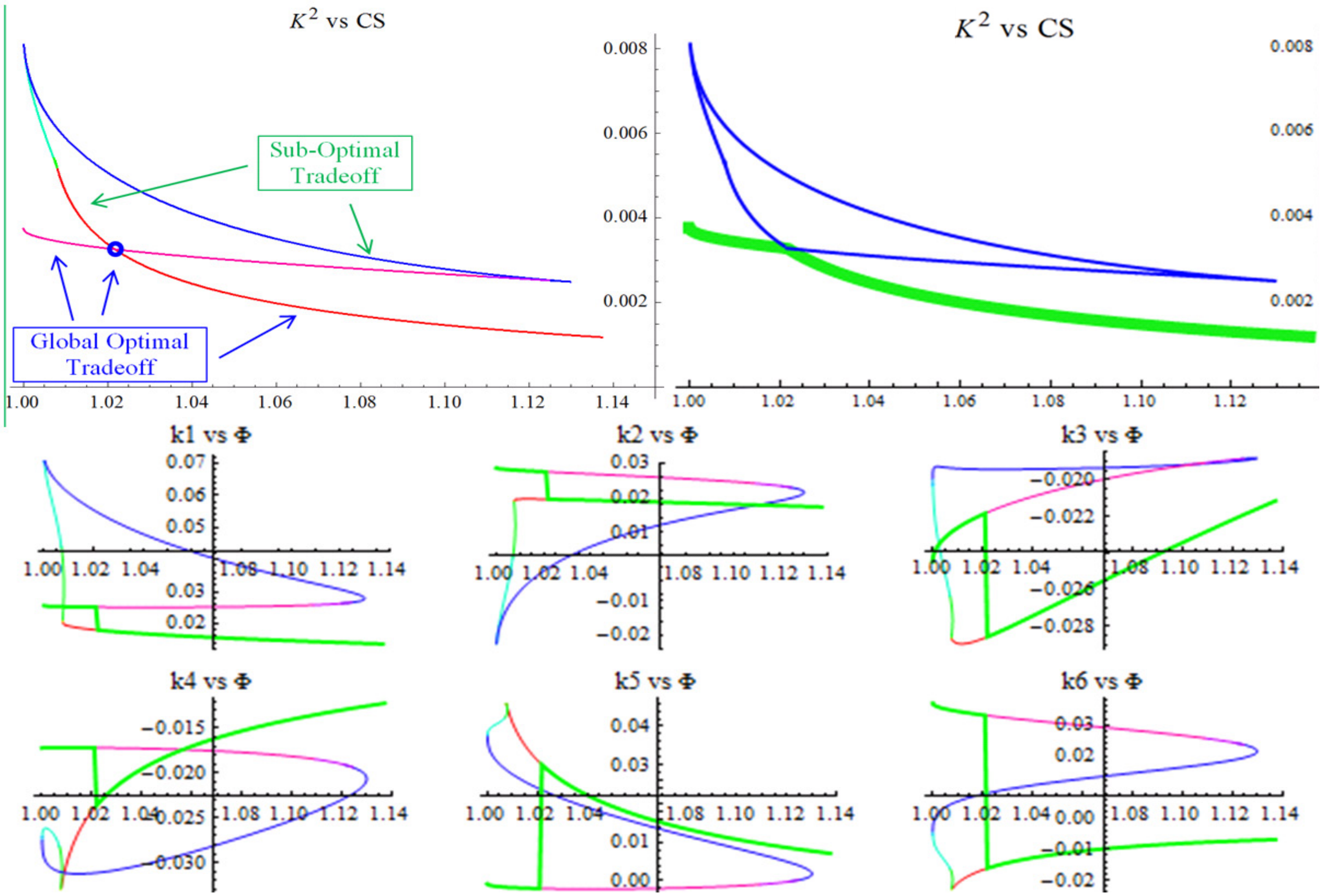

Figure 9. Top Left: Extracting Globally Optimal Tradeoff from solution curve in the $h$–$f$ plane. This is done by joining the red section with the magenta section at the intersection (blue circle), and discard everything above and to the right.

Top Right: Globally Optimal Tradeoff in thick green line. Everything else (blue) is Pareto-dominated by this.

Bottom: The 6 quadrupoles of Figure 8 following the path of globally optimal tradeoff indicated by green highlights. Start with the red curve in Figure 8, then jump to the magenta curve at $f(CS)$=1.02, the blue circle location in the top left plot. An equivalent plot can be made in the $k_m$ vs $h$ ($K^2$) plane as well. The fact that $\lambda$ is negative everywhere makes this process unambiguous for both $f$ and $h$. The symbol Φ is the same as CS here.

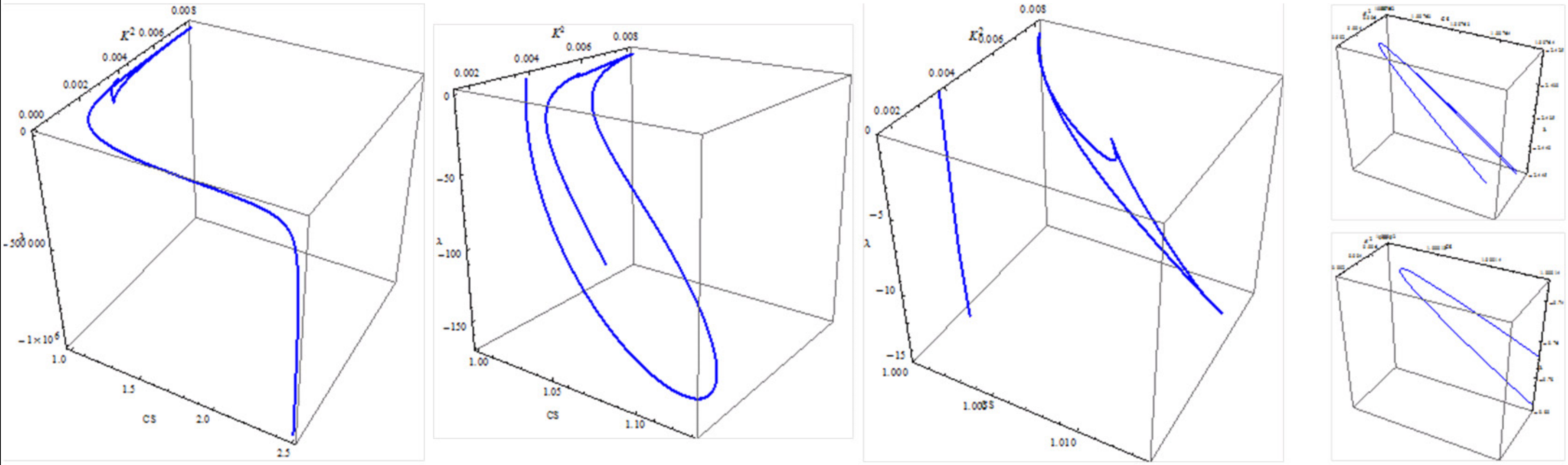

Figure 10. Evolution of the solution curve corresponding to Figures 3 and 6 in the 3D space spanned by $\lambda$–$f$–$h$, zoomed to various subsections.

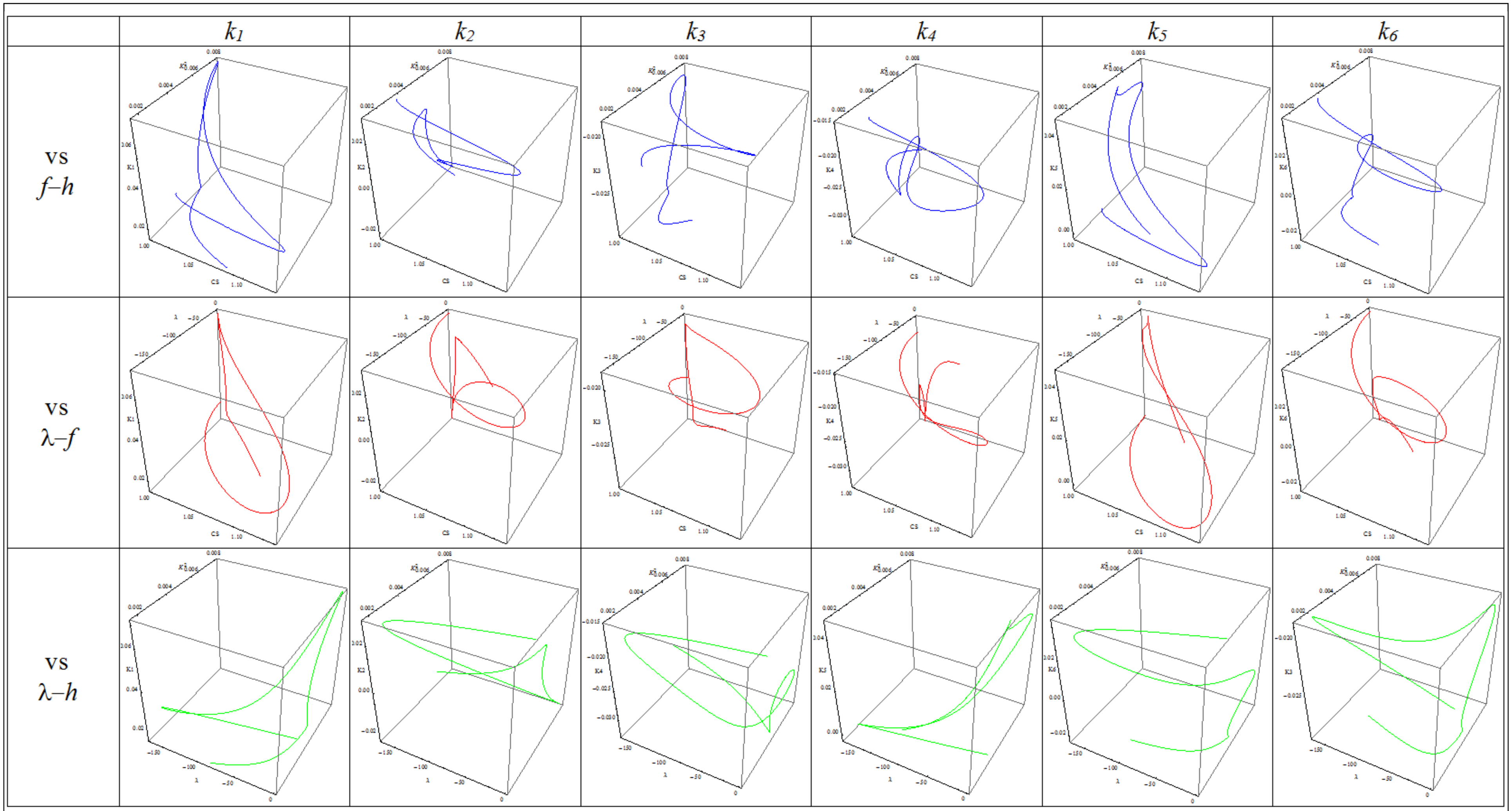

Figure 11. Evolution of the 6 solutions for $k_m$ corresponding to Figures 3 and 6, with respect to two of the 3 parameters λ, $f$, and $h$. Only end sections are shown.

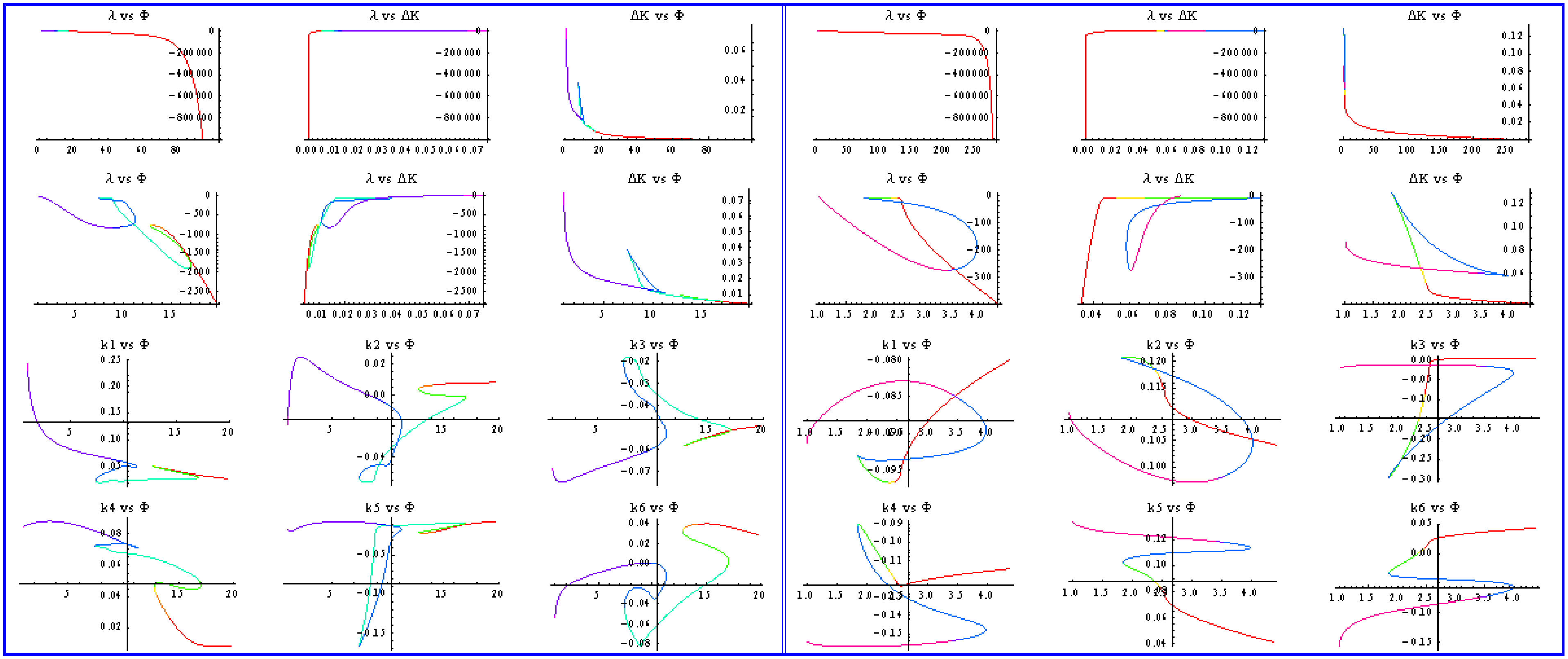

Figure 12. Evolution of λ, $f$, and $h$, as well as all $k_m$'s for path corresponding to constrained optimal solutions to correct optical transport errors.
Left: 30° phase advance optics with introduced error. Optimization is subject to minimal RMS $\delta k_m$ (See eqn. (19)).
Right: 30° phase advance optics with introduced error. Optimization is subject to minimal RMS $k_m$.
Despite the convoluted pattern of the solution paths, inspection of the $f$–$h$ graph (upper right in each case) shows extracting the Pareto front should be straightforward in both cases, thanks to the fact that the slope (=λ) is always negative.

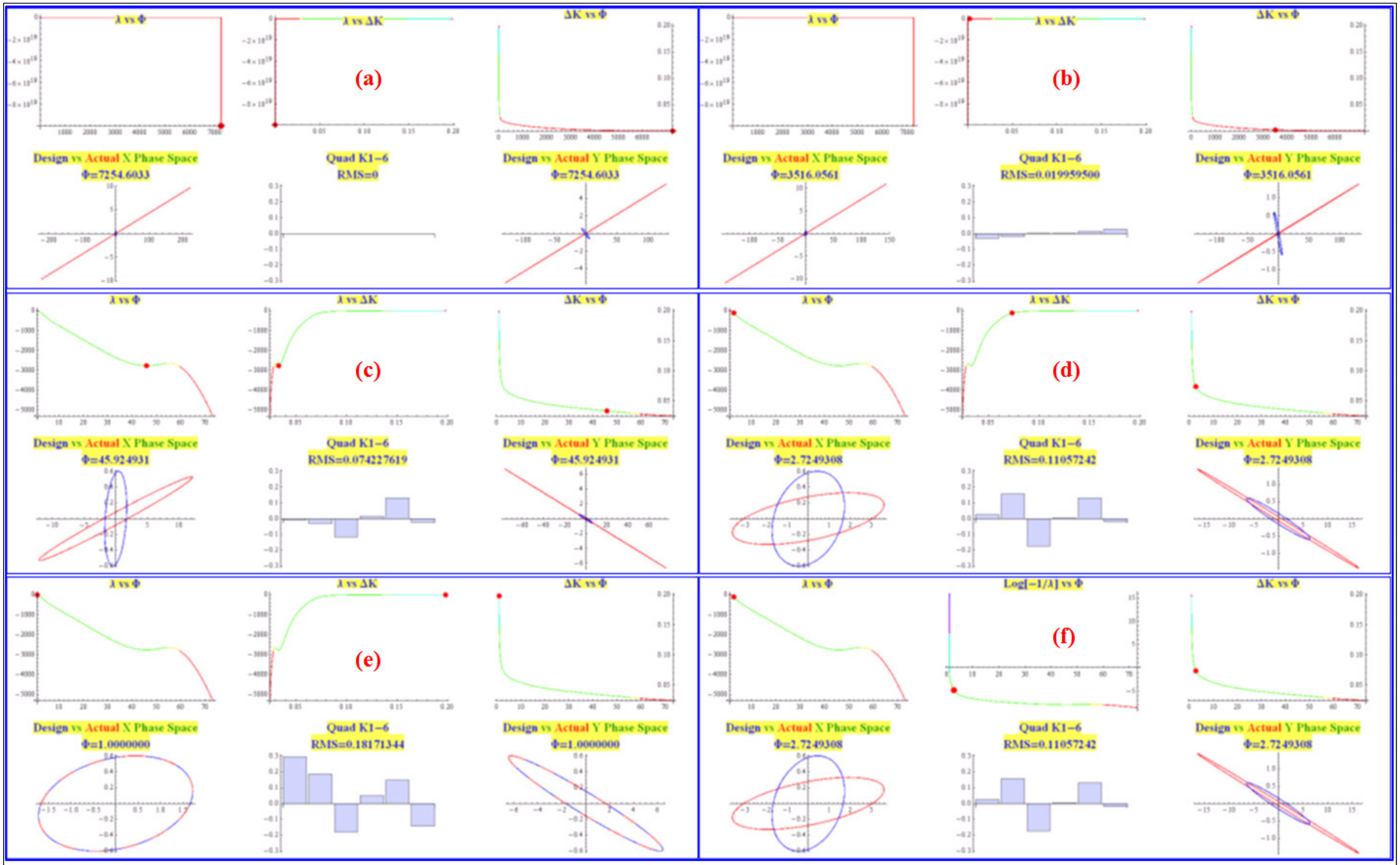

Figure 13. Details of constrained optimization to correct transport errors in a 120° phase advance optics. Each graph represents a particular point (red dot) in the solution path. The ellipses in the bottom plots represent target (blue) and intermediate (red) beam phase space contour at a given correction stage. Complete correction happens when the two ellipses coincide. The bar chart shows the strengths of all 6 quadrupoles at that stage. (a): Initial mismatch factor $\Phi$ reaches 7254! (b): Small amount of quadrupole correction brought $\Phi$ to 3516. (c): Inflection point (local minimum in $\lambda$–$\Phi$ plot), or point of diminished return reached, beyond which large quadrupole change is need for small gain in $\Phi$. (d): Further up the curve past point of diminished return. (e): 100% correction with large quadrupole expense at very last stage. (f): Comparison between $\lambda$ and the *f*–*h* curve.

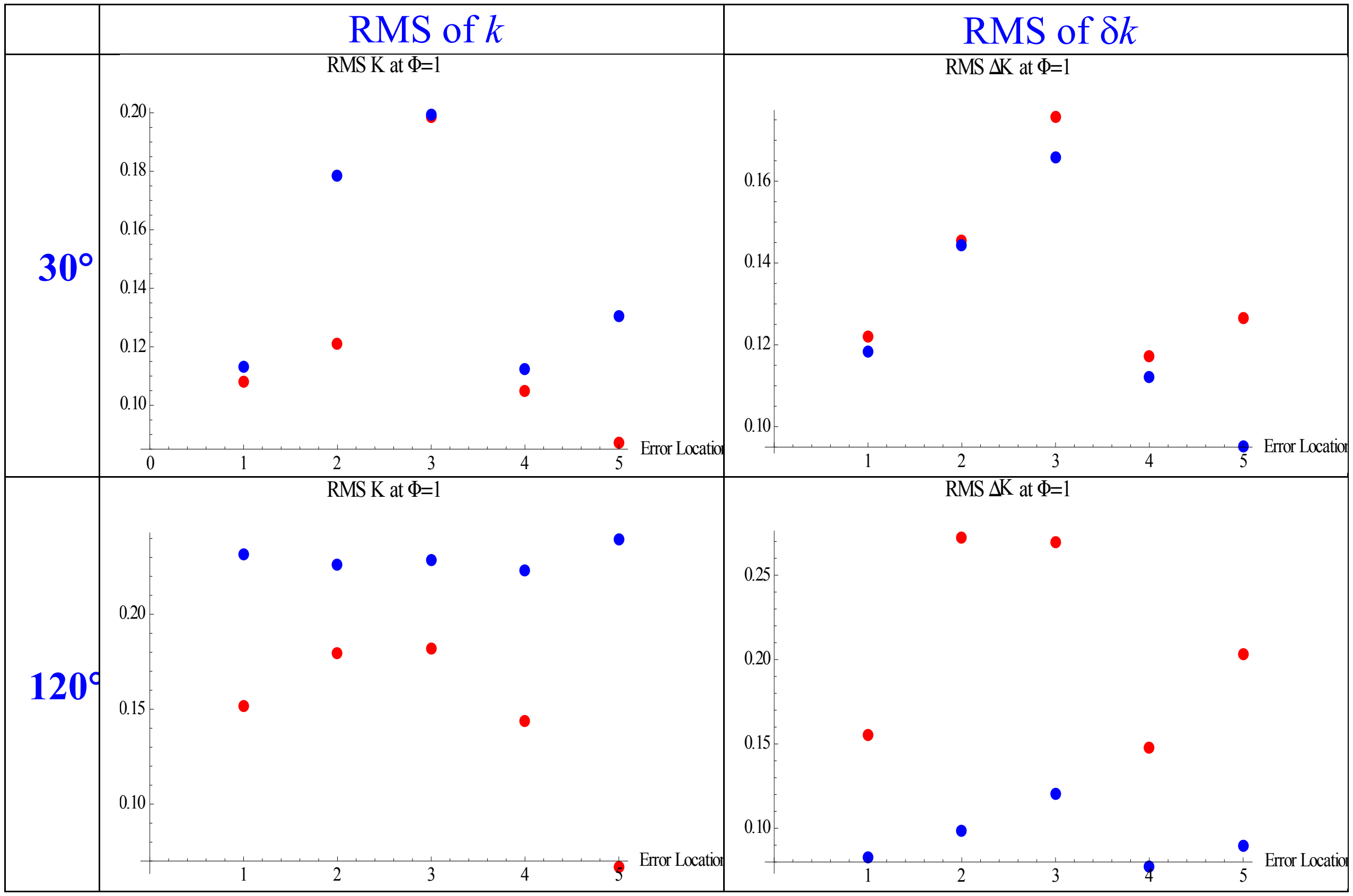

Figure 14. Result of correcting various modes of transport optics error through optimization using proposed algorithm under constraint of either minimal RMS $\delta k_m$ (**blue** dots, see eqn. (19)) or minimal RMS $k_m$ (**red** dots). In all cases the goal of correction (mismatch factor=1, $\lambda$=0) has been achieved, but with varying solutions for $k_m$'s. The RMS values of either final $\delta k_m$ or $k_m$ of the absolute optima (100% correction) are plotted for two different baseline optics (30° & 120°). Apparently the final optimum, although formally constraint independent, depended on the path through which it was arrived at, which in turn depended on the constraint chosen. The final solution reflects the original constraint, namely, if the constraint is on minimal RMS $\delta k_m$, then the final optimal solution will have a smaller RMS $\delta k_m$ than that derived from a constraint of minimal RMS $k_m$, and vice versa.

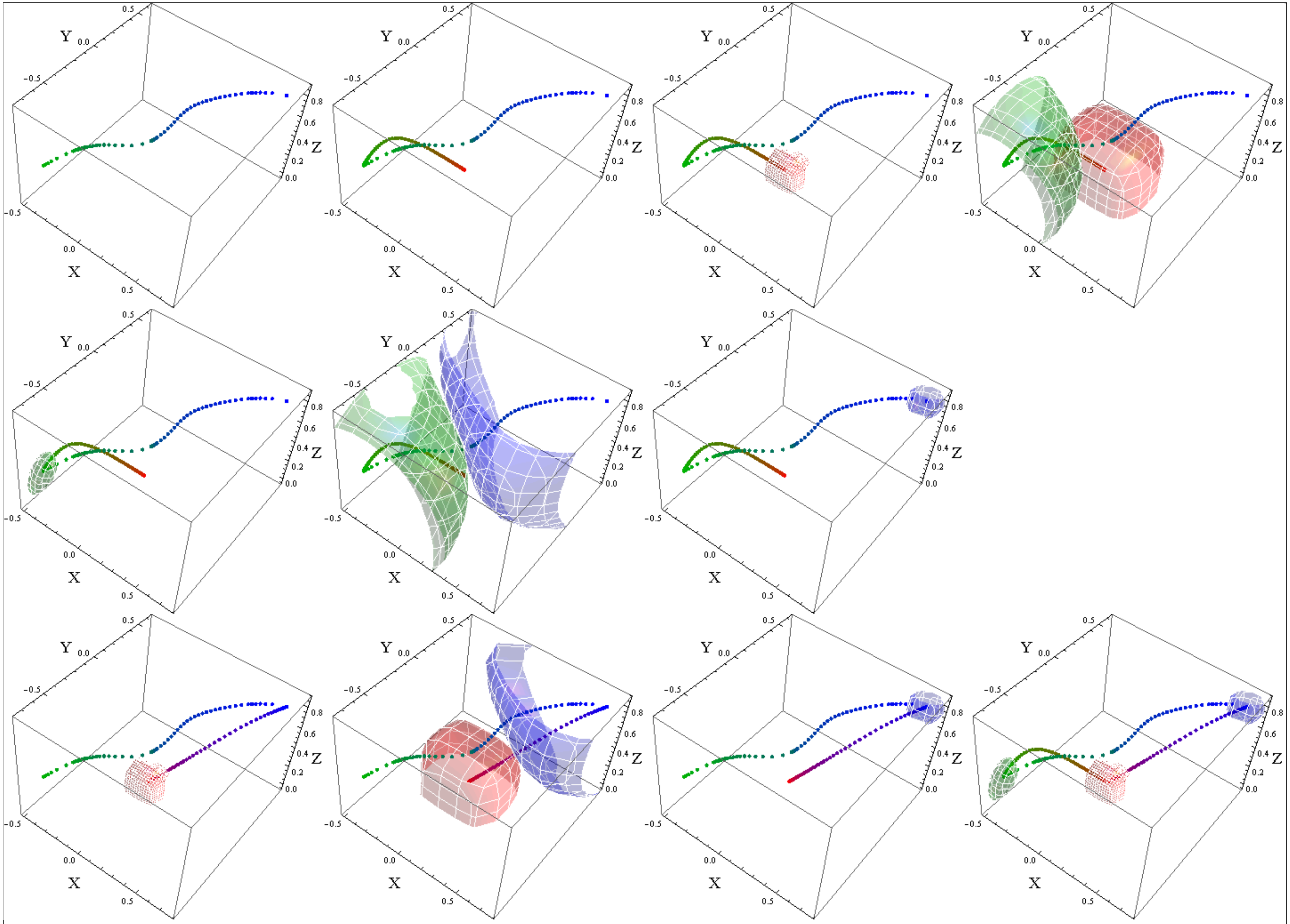

Figure 15. Concept of using artificial constraint to solve for tradeoff curve where neither end point is known a priori (continued on next page).

Figure 15 continued:

A hypothetical problem is depicted here to illustrate the concept of using artificial constraint to solve for tradeoff curve where neither end point is known a priori.

First Row:

The goal is to map out the tradeoff curve between two functions, Green and Blue, with respective CIO's at two ends of the curve shown in the first 3D plot in the variables X, Y & Z. But unlike the quadrupole matching example with trivial constraint, there is no easy way to determine the CIO of either function, thus making it impossible to start from one end and integrate out the rest of the tradeoff curve.

An artificial constraint function, Red, with trivial CIO is introduced ($2^{nd}$ plot), with corresponding tradeoff curve to the CIO of Green. The $3^{rd}$ plot shows the beginning of the integration process with initial equipotential contour of Red, followed by more advanced stage of integration ($4^{th}$ plot).

Second Row:

The previous integration ends up at the CIO of Green (first plot), from which point the integration toward the CIO of Blue can now be launched ($2^{nd}$ plot) until the CIO of Blue is reached ($3^{rd}$ plot). This completes the process of mapping out the tradeoff curve between Green and Blue.

Second Row:

The entire process can equally well start by integrating from Red towards Blue first ($1^{st}$ , $2^{nd}$ & $3^{rd}$ plots), followed by backward integration from Blue to Green instead. All 3 functions can be linked by forward and backward integration along the tradeoff curve ($4^{th}$ plot).

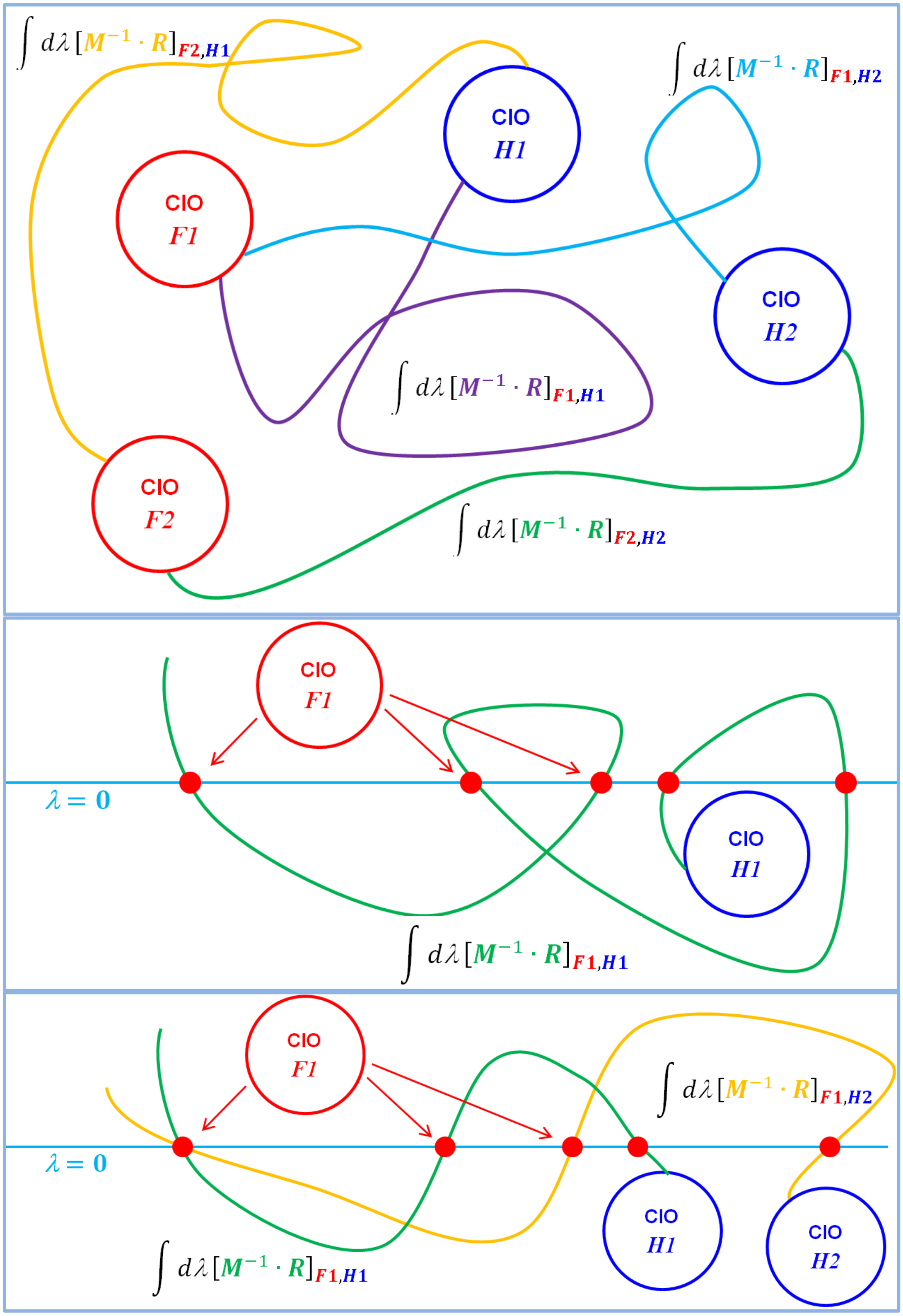

Figure 16. Some possible global structures of the tradeoff curves.

Top: CIO's for different $F$ and $H$ linked by tradeoff curves defined according to (9) and (12). $\int d\lambda\,[\boldsymbol{M}^{-1}\cdot\boldsymbol{R}]_{F,H}$ is symbolic for the tradeoff curve linking $F$ and $H$.

Middle: Tradeoff curve is uniquely fixed for given $F$ and starting CIO of $H$. This curve however can cross $\lambda$=0 multiple times, each time resulting in a distinct CIO for $F$. On the other hand, will the curve launched from one CIO of $H$ visit all CIO's of $F$?

Bottom: It is conceivable, although cases remain to be made, that the potentially multiple CIO's of $F$ can be linked by different tradeoff curves to different CIO's of either the same or different $H$'s. The topological structure of this problem appears highly interesting.

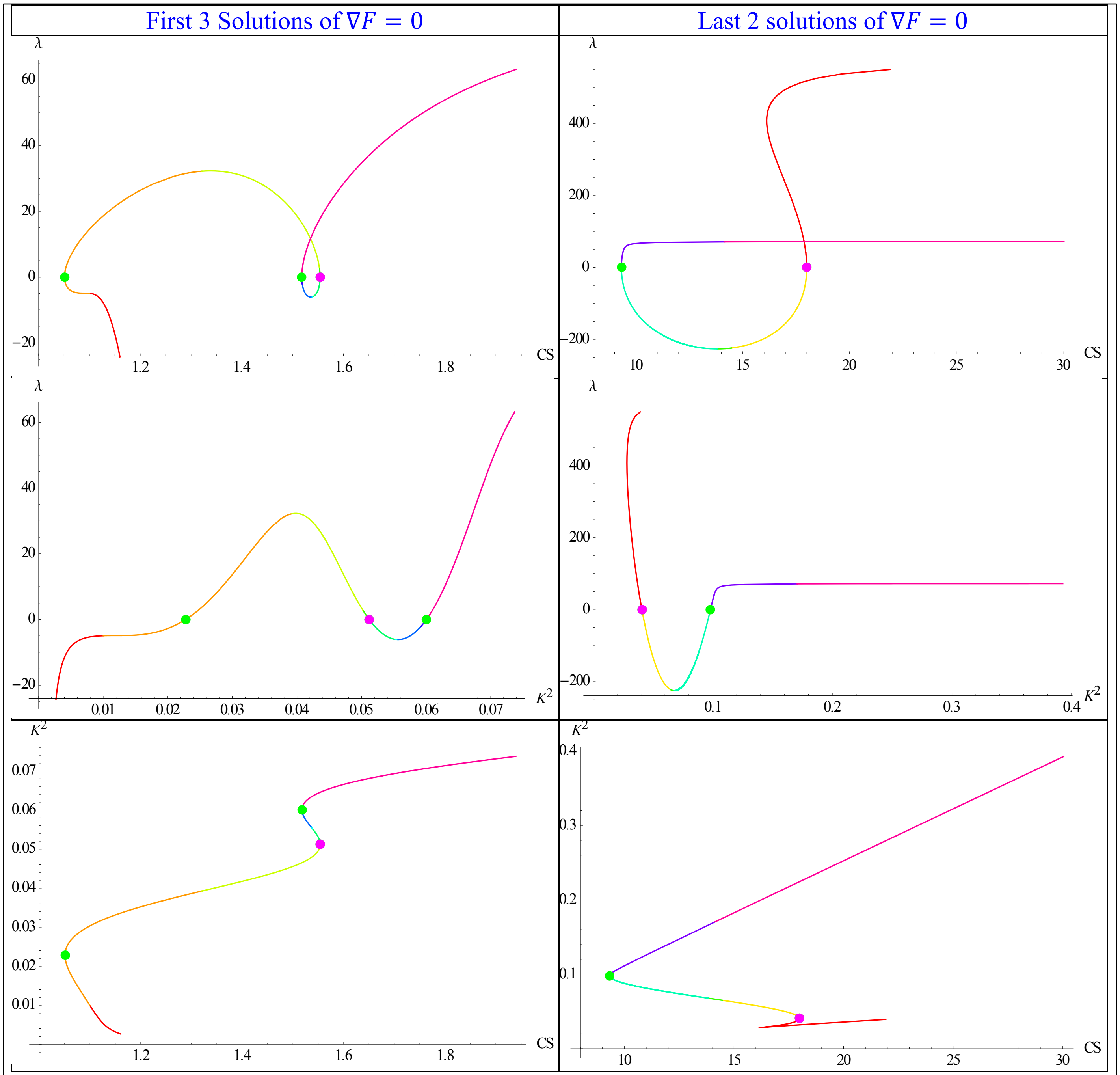

Figure 17. Extending the solution path of tradeoff curve in a realistic betatron matching problem beyond $\lambda = 0$ shows its crossing of more solutions to (5), $\nabla F = 0$. Independent computation identified 5 real solutions to (5), 3 minima and 2 saddle points (see footnote 26). The plots in the above table show that these solution points are visited by two segments of the tradeoff curve. The two segments may also be connected, possibly at infinity. However the disparate scales of these two sets of solutions made establishing such a connection numerically difficult.

Left Column: First 3 solutions to $\nabla F = 0$. Green dots are minima. Magenta dots are saddle points/maxima depending on context (see footnote 26).
Right Column: Last 2 solutions to $\nabla F = 0$.

Top Row: Tradeoff curve ($-\infty$ section not shown) in space of $\lambda$ vs $F$ taken beyond $\lambda = 0$ into further solutions. It passes more solutions whenever it crosses $\lambda = 0$, including both minima and 2 saddle points.
Middle Row: Tradeoff curve in space of $\lambda$ vs $H$ taken beyond $\lambda = 0$.
Bottom Row: Tradeoff curve in space of $H$ vs $F$ taken beyond $\lambda = 0$.

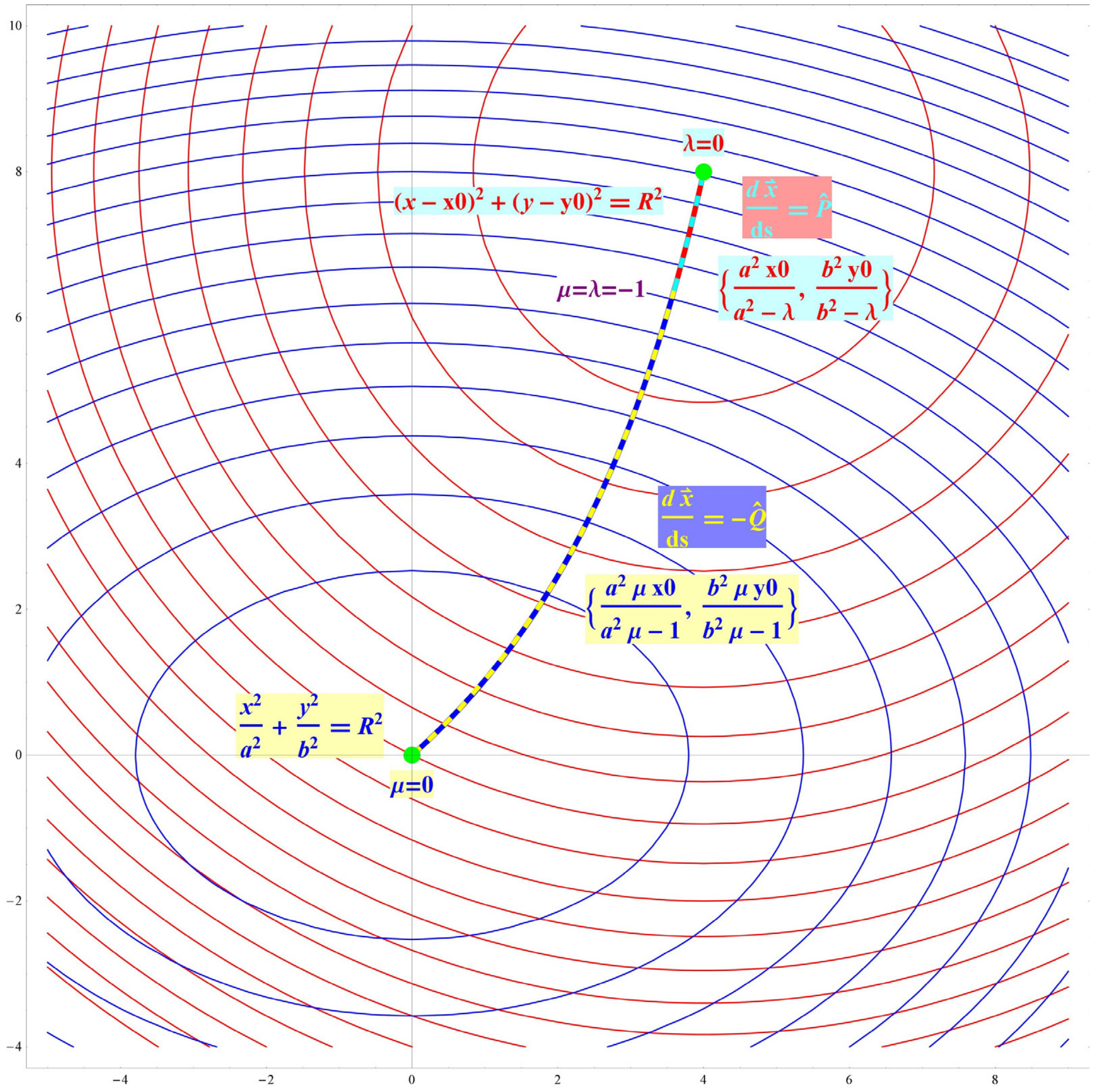

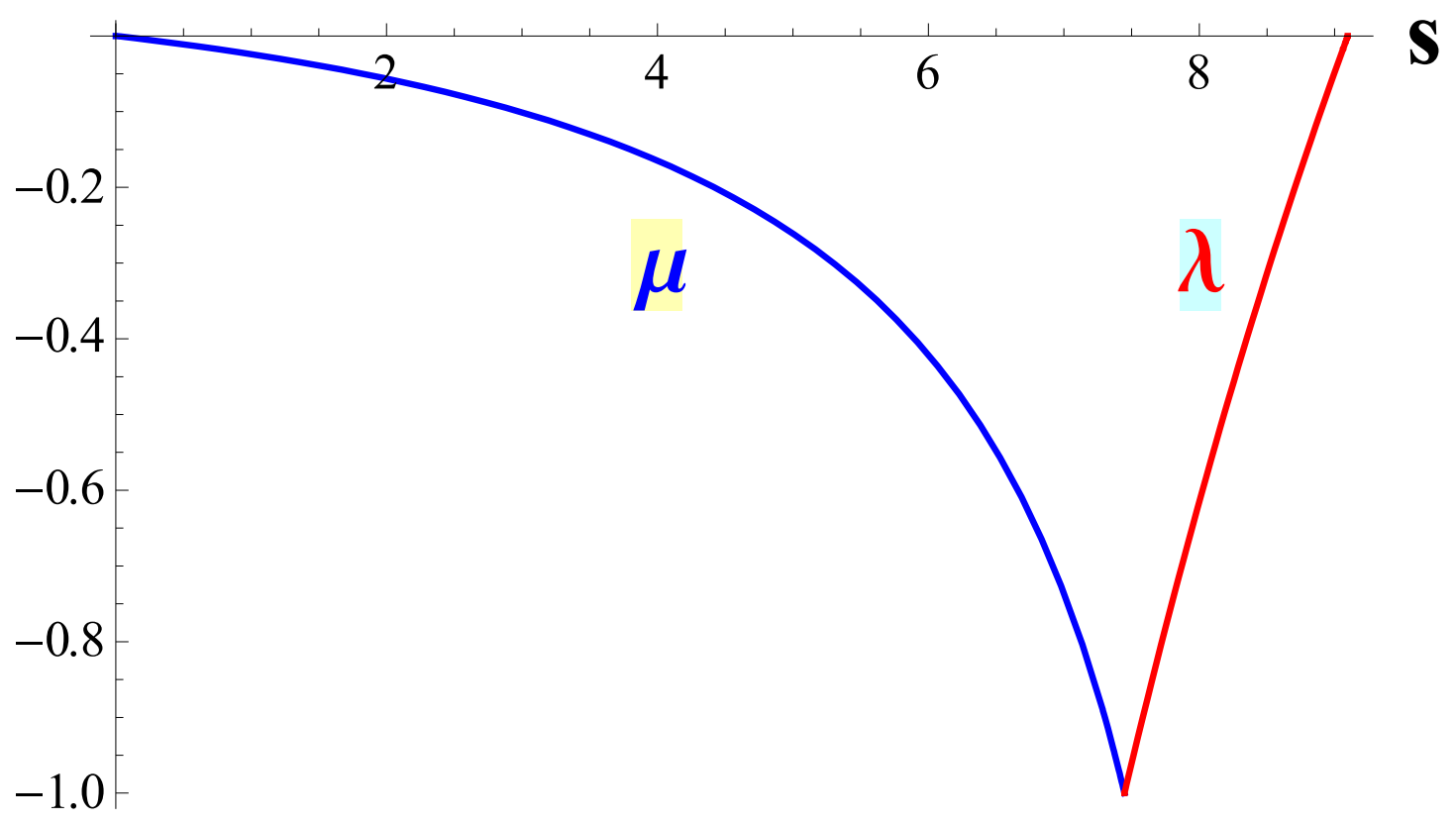

Figure 18. Example showing agreement between analytical expression of the tradeoff curve (yellow & cyan lines) between minima of two competing objectives, and numerically integrated results (blue and red dotted) based on (23), (24) and (25). Evolution of $\mu$ & $\lambda$ are shown below as function of path length along the tradeoff curve.

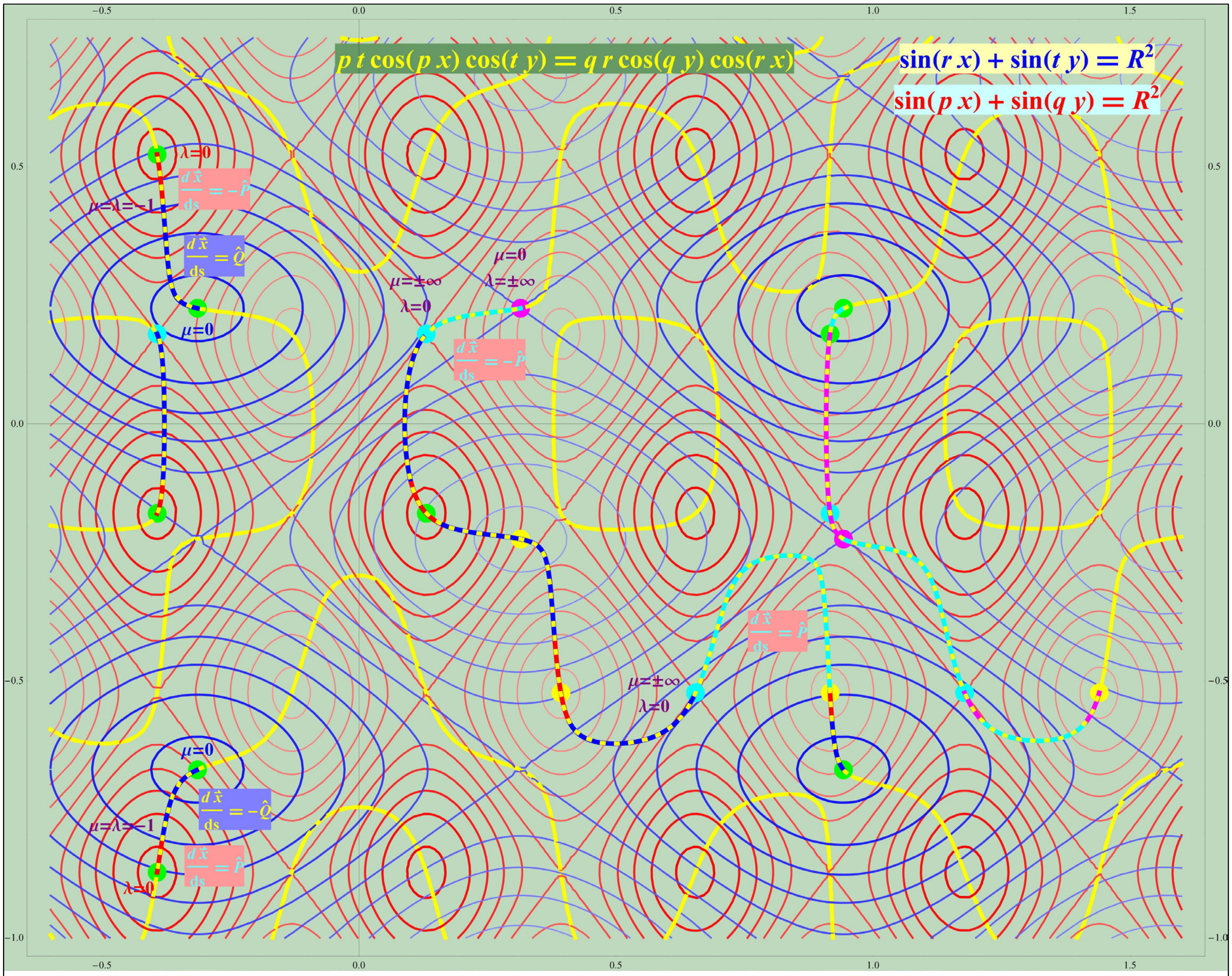

Figure 19. Performance of integration system (23), (24) and (25) in an example where the two competing objectives are given by the blue and red equations displayed at the upper right hand corner, with corresponding color coded equipotential contours. The brightness of the contour lines indicates functional value, from peak (bright) to trough (dimmed). Tradeoff curves are analytically given by the yellow equation on top, and shown as yellow contours. All dotted lines represent paths produced by integrating (23), (24) and (25) under different conditions. Blue and red lines represent integration paths originating or ending at extrema of either function. In some cases before the next extremum is reached the integration must pass through a saddle point, in which case additional $\lambda/\mu$ switching must take place to avoid artificial singularity. These occurrences are indicated by cyan and magenta lines.

Other color coding: Green dot: peaks of both functions; yellow dot: troughs of both functions; cyan dot: saddle points of red function; magenta dot: cyan dot: saddle points of blue function.